\documentclass[11pt]{amsart}
\usepackage[latin1]{inputenc}
\usepackage[OT2,T1]{fontenc}
\usepackage{lmodern}
\usepackage{amsmath}
\usepackage{amssymb}
\usepackage{mathrsfs}
\usepackage[letterpaper,hcentering,hscale=0.8,top=3.25cm,bottom=3.25cm]{geometry}

\usepackage[ps,dvips,arrow,matrix,tips]{xy}

\entrymodifiers={+!!<0pt,\the\fontdimen22\textfont2>}

\SelectTips{cm}{10}

\newcounter{enonce}[subsection]
\let\myoldsection=\section\renewcommand{\section}[1]{\setcounter{enonce}{0}\myoldsection{#1}}
\makeatletter
\def\theenonce{\ifnum\value{subsection}=0{}\thesection\else\thesubsection\fi.\@arabic\c@enonce}
\def\theequation{\ifnum\value{subsection}=0{}\thesection\else\thesubsection\fi.\@arabic\c@equation}
\numberwithin{equation}{subsection}
\makeatother

\makeatletter
\def\myrightarrow{{\setbox\z@\hbox{$\rightarrow$}\dimen0\ht\z@\multiply\dimen0 6\divide\dimen0 10\ht\z@\dimen0\box\z@}}
\def\myleftarrow{{\setbox\z@\hbox{$\leftarrow$}\dimen0\ht\z@\multiply\dimen0 6\divide\dimen0 10\ht\z@\dimen0\box\z@}}
\def\myrightarrowfill@{\arrowfill@\relbar\relbar\myrightarrow}
\def\myleftarrowfill@{\arrowfill@\myleftarrow\relbar\relbar}
\newcommand{\myxrightarrow}[2][]{\ext@arrow 0359\myrightarrowfill@{#1}{#2}}
\newcommand{\myxleftarrow}[2][]{\ext@arrow 3095\myleftarrowfill@{#1}{#2}}
\makeatother
{\setbox0\hbox{$ $}}
\lineskip=\lineskiplimit
\newenvironment{enonce}[1]{\noindent{\textbf{#1}}---\!\, \begin{itshape}}{\end{itshape}}
\newenvironment{proposition}[1][]{\refstepcounter{enonce}\begin{enonce}{Proposition \theenonce{}#1 }}{\end{enonce}}
\newenvironment{theorem}[1][]{\refstepcounter{enonce}\begin{enonce}{Theorem \theenonce{}#1 }}{\end{enonce}}
\newenvironment{corollary}[1][]{\refstepcounter{enonce}\begin{enonce}{Corollary \theenonce{}#1 }}{\end{enonce}}
\newenvironment{lemma}[1][]{\refstepcounter{enonce}\begin{enonce}{Lemma \theenonce{}#1 }}{\end{enonce}}
\newenvironment{construction}[1][]{\refstepcounter{enonce}\noindent{\textbf{Construction \theenonce{}#1 }}---\!\, }{}
\newenvironment{remark}[1][]{\refstepcounter{enonce}\noindent{\textbf{Remark \theenonce{}#1 }}---\!\, }{}
\newenvironment{question}[1][]{\refstepcounter{enonce}\noindent{\textbf{Question \theenonce{}#1 }}---\!\, }{}
\renewenvironment{proof}[1][]{\noindent{\it{}Proof{}#1 }--- }{\hfill $\square$}
\newenvironment{sketchofproof}[1][]{\noindent{\it{}Sketch of proof{}#1 }--- }{\hfill $\square$}

\input cyracc.def
\DeclareFontFamily{U}{russian}{}
\DeclareFontShape{U}{russian}{m}{n}
        { <5><6> wncyr5
        <7><8><9> wncyr7
        <10><10.95><12><14.4><17.28><20.74><24.88> wncyr10 }{}
\DeclareSymbolFont{Russian}{U}{russian}{m}{n}
\DeclareSymbolFontAlphabet{\mathcyr}{Russian}
\makeatletter
\let\@math@cyr\mathcyr
\renewcommand{\mathcyr}[1]{\@math@cyr{\cyracc #1}}
\makeatother

\newcommand{\Ba}{{\mathcyr{B}}}

\newcommand{\Sm}{\mathrm{Sm}}
\renewcommand{\phi}{\varphi}
\renewcommand{\epsilon}{\varepsilon}
\renewcommand{\emptyset}{\varnothing}
\newcommand{\limind}{\varinjlim}
\newcommand{\limproj}{\varprojlim}

\newcommand{\Hom}{{\mathrm{Hom}}}
\newcommand{\Ext}{{\mathrm{Ext}}}

\newcommand{\catExt}{{\mathbf{{A}}}}
\newcommand{\Br}{{\mathrm{Br}}}
\newcommand{\Arond}{\mathscr{A}}
\newcommand{\Lrond}{\mathscr{L}}
\newcommand{\Orond}{\mathscr{O}}
\newcommand{\Prond}{\mathscr{P}}
\newcommand{\Qrond}{\mathscr{Q}}
\newcommand{\Urond}{\mathscr{U}}
\newcommand{\C}{\mathbf{C}}
\newcommand{\Q}{\mathbf{Q}}
\newcommand{\Z}{\mathbf{Z}}
\newcommand{\N}{\mathbf{N}}
\renewcommand{\P}{\mathbf{P}}
\newcommand{\A}{\mathbf{A}}
\newcommand{\ensemble}[2]{\{ #1 \, ; \, #2 \}}
\newcommand{\uplet}[2]{#1, \mskip2.5mu \ldots \mskip-1mu, \mskip2.5mu #2}
\newcommand{\Gal}{\mathrm{Gal}}
\newcommand{\Alb}{\mathrm{Alb}}
\newcommand{\TAlb}{\mathrm{TAlb}}
\newcommand{\T}{\mathrm{T}}
\newcommand{\Div}{\mathrm{Div}}
\newcommand{\Pic}{\mathrm{Pic}}
\newcommand{\Rstar}{\mathrm{R}^\star}
\newcommand{\End}{\mathrm{End}}
\newcommand{\Id}{\mathrm{Id}}
\newcommand{\Sh}{\mathrm{Sh}}
\newcommand{\NS}{\mathrm{NS}}
\newcommand{\Gm}{\mathbf{G}_\mathrm{m}}
\newcommand{\GmS}{\mathbf{G}_{\mathrm{m},S}}
\newcommand{\GmA}{\mathbf{G}_{\mathrm{m},A}}
\newcommand{\longisoto}{\myxrightarrow{\,\,\sim\,\,}}
\newcommand{\isoto}{\myxrightarrow{\,\sim\,}}
\newcommand{\kbar}{{\mkern1mu\overline{\mkern-1mu{}k\mkern-1mu}\mkern1mu}}
\newcommand{\Kbar}{{\mkern1mu\overline{\mkern-1mu{}K\mkern-1mu}\mkern1mu}}
\newcommand{\kalg}{k^{\mathrm{alg}}}
\newcommand{\pitilde}{{\widetilde{\pi}}}
\newcommand{\alphatilde}{{\widetilde{\alpha}}}
\DeclareMathOperator{\Spec}{Spec}
\DeclareMathOperator{\Coker}{Coker}
\newcommand{\ob}{{{ob}}}

\newcommand{\EM}{{{E}^{\mathcal{M}}}}
\newcommand{\E}{{{E}}}
\newcommand{\R}{{\mathrm{R}}}
\newcommand{\Pgen}{P_{\mathrm{gen}}}

\author{Olivier Wittenberg}
\title[On Albanese torsors and the elementary obstruction]{On Albanese torsors and the elementary obstruction}
\date{September 26, 2006; revised on August 15, 2007}
\address{Department of Mathematics, Rice University, 6100 S. Main St., Houston, TX 77251-1892, USA}
\email{olivier.wittenberg@rice.edu}

\begin{document}
\maketitle

\section{Introduction}

If~$X$ is a principal homogeneous space, or torsor, under a semi-abelian
variety~$A$ over a field~$k$, one can naturally associate with it two integers
which give a measure of its nontriviality:
the index~$I(X)$,
defined as the greatest common divisor of the degrees of the closed points
of~$X$, and the period~$P(X)$, which is the order of the class of~$X$ in the
Galois cohomology group $H^1(k,A)$.  The period divides the index.
In the classical situation where~$A$ is an abelian variety, the extent to which the period and the index can fail to be equal has been investigated by
a number of authors, including
Lang and Tate (see~\cite{langtate}), Ogg, Cassels, Shafarevich, Lichtenbaum.
On the other hand, the well-known period-index problem for central simple algebras, first raised by Brauer, is a
particular case of the study of periods and indices of torsors under tori.

Let~$X$ now denote an arbitrary smooth variety over~$k$ (not necessarily complete).   One may define the index~$I(X)$ as above,
and the period~$P(X)$ as the supremum of~$P(Y)$ when~$Y$ ranges over all torsors
under semi-abelian varieties such that there exists a morphism $X \rightarrow Y$.
As was proven by Serre, there is a universal object among
morphisms from~$X$ to torsors under semi-abelian varieties: the Albanese
torsor $X \rightarrow \Alb^1_{X/k}$.  Hence $P(X)=P(\Alb^1_{X/k})$.
Periods of open varieties turn out to be relevant even if one is primarily
interested in complete varieties.  Indeed the period of a small enough dense open
subset~$U$ of a complete variety~$X$ often provides valuable information which
is not given by the period of~$X$ alone.  For instance, if~$X$ is the Severi-Brauer
variety attached to a central simple algebra~$\alpha$, then obviously $P(X)=1$ (since $\Alb^1_{X/k}$ is just a point), whereas
$P(U)$ can be shown to equal the exponent of~$\alpha$ if
$U \subseteq X$ is small enough.

Quite generally, if $U \subseteq X$ is open dense, then~$P(U)$ divides~$I(X)$.  Thus, if we define the generic period~$\Pgen(X)$
as the supremum of~$P(U)$ when~$U$ ranges over all dense open subsets of~$X$, then
for a $0$\nobreakdash-cycle
of degree~$1$ to exist on~$X$, it is necessary that $\Pgen(X)=1$.
Colliot-Thélène and Sansuc defined and studied in~\cite{ctsandesc2}
another general obstruction to the existence of $0$\nobreakdash-cycles of degree~$1$, the so-called elementary obstruction.
Namely, denoting by~$\kbar$ a separable closure of~$k$ and by~$\kbar(X)$ the function field of $X \otimes_k \kbar$,
the elementary obstruction is said to vanish if and only if the inclusion $\kbar^\star \subset \kbar(X)^\star$ admits
a Galois-equivariant retraction.

The theme of the present paper is that these two obstructions are very tightly related.  We shall prove in particular
that $\Pgen(X)=1$ implies the vanishing of the elementary obstruction, over an arbitrary field~$k$ (Theorem~\ref{thorder}).
When~$k$ is a $p$\nobreakdash-adic field, a real closed field, a number field, or a field of dimension~$\leq 1$ in the sense of Serre, we are able to obtain
stronger statements.  Most notably, the condition $\Pgen(X)=1$ is even equivalent to the vanishing of the elementary obstruction
if~$k$ is a $p$\nobreakdash-adic or real closed field (Theorem~\ref{thapploc}), or if~$k$ is a number field and the Tate-Shafarevich group of the Picard variety of~$X$ is finite
(Theorem~\ref{mainapp}); we give several applications of these results in §\ref{sectarith}.
On the other hand, if~$k$ is a field of dimension~$\leq 1$, we prove that the elementary obstruction always vanishes (Theorem~\ref{thcdone}), while it is well known
that~$P(X)$ (and therefore~$\Pgen(X)$) may be greater than~$1$.

The elementary obstruction to the existence of $0$\nobreakdash-cycles of degree~$1$ on~$X$
is known to vanish if and only if the Yoneda equivalence class of a certain $2$\nobreakdash-extension~$e(X)$ of~$\Pic(X \otimes_k \kbar)$ by~$\kbar^\star$
in the category of discrete Galois modules is trivial (see \cite[Proposition~2.2.4]{ctsandesc2}).
One is naturally led to consider an analogous $2$\nobreakdash-extension~$E(X)$ of the relative Picard functor~$\Pic_{X/k}$ by the multiplicative group~$\Gm$;
it is a straightforward generalisation of the elementary obstruction.
We devote the last section of this paper to showing that the $2$\nobreakdash-extension~$E(X)$ contains enough information to reconstruct, in a strong sense,
the Albanese torsor~$\Alb^1_{U/k}$ for any open $U \subseteq X$ (Theorem~\ref{themaintheorem}).
More precisely, for any open $U \subseteq X$, the $2$\nobreakdash-extension~$E(X)$ gives rise, by pullback, to a $2$\nobreakdash-extension~$\EM(U)$ of
the Picard $1$\nobreakdash-motive of~$U$ by~$\Gm$.
Theorem~\ref{themaintheorem} essentially states that~$\Alb^1_{U/k}$ is the variety which parametrises, in an appropriate sense,
``Yoneda trivialisations'' of the $2$\nobreakdash-extension~$\EM(U)$.  In particular, the Yoneda equivalence class of the
$2$\nobreakdash-extension~$\EM(U)$ is trivial if and only if $\Alb^1_{U/k}(k)\neq\emptyset$,
that is, if and only if~$P(U)=1$.
From this we deduce, without any assumption on the field~$k$, that the Yoneda equivalence class of~$E(X)$ is trivial if and only if $\Pgen(X)=1$
(Corollary~\ref{quatcor}).  This is of course stronger than our previous assertion that the elementary obstruction vanishes if~$\Pgen(X)=1$.

Producing trivialisations of the $2$\nobreakdash-extension~$\EM(U)$ from rational points of~$\Alb^1_{U/k}$ (or from rational points of~$U$)
is not particularly difficult (see Proposition~\ref{propyondeq} and the proof of Corollary~\ref{quatcor});
it is however much less clear how to convert trivialisations of~$\EM(U)$ into rational points of~$\Alb^1_{U/k}$.
One of the key tools in the proof of Theorem~\ref{themaintheorem}
is a result which might be of some independent interest: we exhibit an explicit Poincaré sheaf
on any abelian variety (Theorem~\ref{thpoincare}).

\bigskip
\emph{Acknowledgements.}
This work grew out of an attempt to answer the questions raised by Borovoi,
Colliot-Thélène, and Skorobogatov in~\cite{boctsko},
and to get a better understanding of a result of Skorobogatov~\cite[Proposition~2.1]{skoronote}.
In response to my queries about automorphisms of $n$\nobreakdash-extensions,
Shoham Shamir very kindly directed me to the papers of Retakh~\cite{retakh} and of Neeman and Retakh~\cite{neemanretakh},
and explained to me a variant of Lemma~\ref{lemmashoham} below.
I am grateful to Dennis Eriksson for our many stimulating discussions about Albanese torsors and generic periods;
it should be noted that generic periods first appeared in~\cite{dennis}.
Finally I~am pleased to thank Jean-Louis Colliot-Thélène for his constant interest and encouragement,
and Joost van Hamel for a number of useful comments and for pointing out the connections between §\ref{secquat} and~\cite{vanhamel}.

\bigskip
\emph{Notation.}
Let~$k$ be a field. A variety over~$k$ is by definition a $k$\nobreakdash-scheme of finite type.
The Néron-Severi group~$\NS(X)$ of a smooth $k$\nobreakdash-variety~$X$ is the quotient of~$\Pic(X)$ by the subgroup of all
classes which are algebraically equivalent to~$0$ over~$k$.
If~$A^0$ is an algebraic group over~$k$ and~$A^1$ is a $k$\nobreakdash-torsor under~$A^0$, we shall generally denote by~$A^n$
the $n$\nobreakdash-fold contracted product of~$A^1$ with itself under~$A^0$ (so that the class of~$A^n$ in the Galois cohomology group $H^1(k,A^0)$
is~$n$ times the class of~$A^1$).  If~$A$ is a semi-abelian variety over~$k$ (that is, an extension of an abelian variety by a torus), we denote by~$\T{}A$
the largest $k$\nobreakdash-torus contained in~$A$ (so that $A/\T{}A$ is an abelian variety).
Let~$\Arond$ be an abelian category. We denote by~$C(\Arond)$ the category of complexes of objects of~$\Arond$.
If~$f$ is a morphism in~$C(\Arond)$, we also denote by~$C(f)$ the mapping cone of~$f$ (\cite[§1.4]{kascha}); the context will make it clear which
meaning is intended.
For $n \in \N$ and $A,B \in \Arond$ (resp.~$A,B\in C(\Arond)$), we denote by
$\Ext^n_\Arond(A,B)$
the corresponding $\Ext$ (resp.~hyperext) group.
If~$\Arond$ is the category of discrete modules over a ring~$R$ (resp.~over a profinite group~$\Gamma$), we
write~$\Ext^n_R(A,B)$ (resp.~$\Ext^n_\Gamma(A,B)$) for~$\Ext^n_\Arond(A,B)$.
Finally, by a $1$\nobreakdash-motive we shall always mean a Deligne $1$\nobreakdash-motive, \emph{i.e.}, a $1$\nobreakdash-motive as defined in~\cite{raynaudmonodr}.

\section{The order of the elementary obstruction}
\label{sectiontheorder}

Let~$k$ be a field and~$X$ be a geometrically integral variety over~$k$.

The \emph{Albanese variety} and the \emph{Albanese torsor} of~$X$ over~$k$, if
they exist, are a semi-abelian variety $\Alb^0_{X/k}$ over~$k$
and a $k$\nobreakdash-torsor $\Alb^1_{X/k}$ under $\Alb^0_{X/k}$, endowed
with a $k$\nobreakdash-morphism $u_X \colon X \rightarrow \Alb^1_{X/k}$.  They are characterised
by the following universal property: for any semi-abelian variety~$A^0$ over~$k$,
any torsor~$A^1$ under~$A^0$, and any $k$\nobreakdash-morphism $m \colon X \rightarrow A^1$,
there exists a unique $k$\nobreakdash-morphism of varieties
$f^1 \colon \Alb^1_{X/k} \rightarrow A^1$
such that $f^1 \circ u_X = m$, and there exists a unique $k$\nobreakdash-morphism
of algebraic groups $f^0 \colon \Alb^0_{X/k} \rightarrow A^0$
such that $f^1$ is $f^0$\nobreakdash-equivariant.

Obviously, the Albanese variety, the Albanese torsor, and the morphism~$u_X$ are unique up to a unique isomorphism.
Their existence was proven by Serre~\cite{serremunivalb} in case~$k$ is algebraically closed.
By Galois descent it follows that they exist as soon as~$k$ is perfect (see \cite[Theorem~2.1]{dennis} and~\cite[Ch.~V, §4]{serregalg}).
As explained in the appendix to this paper, Serre's arguments can be made to work over separably closed
fields (Theorem~\ref{thalbexiste}).  Galois descent then implies that the Albanese variety,  the Albanese torsor, and the morphism~$u_X$
always exist, without any assumption on~$k$.

We henceforth assume that~$X$ is smooth over~$k$.

If~$X$ is proper over~$k$, the Albanese variety just defined coincides with the classical Albanese variety.
In other words, the semi-abelian variety $\Alb^0_{X/k}$ is then an abelian variety, and as such, it is canonically dual
to the Picard variety $\Pic^0_{X/k,\mathrm{red}}$ of~$X$ over~$k$ (see \cite[Théorème~3.3]{tdte6}; note that $\Pic^0_{X/k,\mathrm{red}}$ is an abelian
variety by \cite[Corollaire~3.2]{tdte6}).

The \emph{period} of~$X$ over~$k$, denoted $P(X)$, is by definition the order of
the class of the torsor $\Alb^1_{X/k}$
in the Galois cohomology group $H^1(k,\Alb^0_{X/k})$.  The \emph{generic period} $\Pgen(X)$ of~$X$ over~$k$ is the supremum of~$P(U)$ when~$U$ ranges
over all dense open subsets of~$X$.
The \emph{index} $I(X)$ of~$X$ over~$k$ is the lowest positive degree of a $0$\nobreakdash-cycle on~$X$.

\bigskip
\begin{proposition}
\label{propdivisibility}
The period, the generic period, and the index satisfy the divisibility relations
$$P(X) \;|\; \Pgen(X) \;|\;
I(X)\rlap{\text{.}}$$
\end{proposition}

\begin{proof}
It suffices to prove that~$P(X)$ divides~$I(X)$.  Indeed, by applying this result to all dense open subsets $U \subseteq X$, one deduces
that $\Pgen(X)$ divides~$I(X)$, thanks to the well-known fact that
$I(U)=I(X)$ for any dense open subset $U \subseteq X$ (see \cite[p.~599]{colliotfinitude}).
For~$P(X)$ to divide~$I(X)$, it suffices that $\Alb^n_{X/k}(k) \neq \emptyset$ for each $n>0$ such that~$X$ contains a closed point of degree~$n$,
where~$\Alb^n_{X/k}$ denotes the $n$\nobreakdash-fold contracted product of the torsor~$\Alb^1_{X/k}$.
On composing $(u_X)^n \colon X^n \rightarrow (\Alb^1_{X/k})^n$ with the projection $(\Alb^1_{X/k})^n \rightarrow \Alb^n_{X/k}$,
we obtain a morphism $X^n \rightarrow \Alb^n_{X/k}$ which is symmetric, and which therefore factors
through the $n$\nobreakdash-fold symmetric power~$X^{(n)}$ of~$X$ (\cite[Proposition~3.1]{milnejac}).
It now only remains to prove that a closed point of degree~$n$ on~$X$ gives rise to a rational point of~$X^{(n)}$.
This is clear if~$k$ is perfect; in general it follows from the existence
of the Grothendieck-Deligne norm map from the Hilbert scheme of points of degree~$n$ on~$X$ to~$X^{(n)}$ (see~\cite[p.~184]{delignecohsup}).
\end{proof}

\bigskip
Let~$\kbar$ denote a separable closure of~$k$ and let $\Gamma=\Gal(\kbar/k)$.
The \emph{elementary obstruction (to the existence of a $0$\nobreakdash-cycle of degree~$1$ on~$X$ over~$k$)} is
the class $\ob(X) \in \Ext^1_\Gamma(\kbar(X)^\star/\kbar^\star,\kbar^\star)$ of the exact sequence of discrete $\Gamma$-modules
$$
\xymatrix{
0 \ar[r] & \kbar^\star \ar[r] & \kbar(X)^\star \ar[r] & \kbar(X)^\star/\kbar^\star \ar[r] & 0\rlap{\text{,}}
}
$$
where $\kbar(X)$ denotes the function field of $X \otimes_k \kbar$.  It is an easy consequence of Hilbert's Theorem~90 that
the existence of a $0$\nobreakdash-cycle of degree~$1$ on~$X$ forces~$\ob(X)$
to vanish (see \cite[Proposition~2.2.2]{ctsandesc2}).
More generally, the same argument shows that the order of $\ob(X)$ divides~$I(X)$.

Our main result in this section is the following:

\bigskip
\begin{theorem}
\label{thorder}
Let~$X$ be a smooth proper geometrically integral variety over a field~$k$.
The order of $\ob(X)$ divides $\Pgen(X)$.
\end{theorem}

\bigskip
In particular, the elementary obstruction to the existence of a $0$\nobreakdash-cycle of degree~$1$ on~$X$
is even an obstruction to the generic period of~$X$ being equal to~$1$.

The integers $P(X)$, $\Pgen(X)$, $I(X)$, and the order of~$\ob(X)$ satisfy no further divisibility relations than those
given by Proposition~\ref{propdivisibility} and Theorem~\ref{thorder}, as we now show with a few examples.
First of all, the order of~$\ob(X)$ does not divide~$P(X)$ in general: any curve of genus~$0$ without rational points has period~$1$
but $\ob(X)\neq 0$.
Conversely, the order of~$\ob(X)$ is not necessarily divisible by~$P(X)$: if~$X$ is a curve of genus~$1$ without rational points over a field
of dimension~$\leq 1$ (for instance the plane cubic $x^3+ty^3+t^2z^3=0$ over~$\C((t))$), then clearly $P(X)>1$, but Theorem~\ref{thcdone} below implies that $\ob(X)=0$.
The period and the generic period differ in the case of a curve of genus~$0$ without rational points.
Examples where the generic period and the index do not coincide are given by Severi-Brauer varieties attached to central simple algebras for which the exponent and the index are not equal.
One can even give examples of varieties which satisfy $\Pgen(X)=1$ and $I(X)>1$.
According to Theorem~\ref{thapploc} below, it suffices to exhibit a geometrically integral variety~$X$, over a $p$\nobreakdash-adic field,
such that $\ob(X)=0$ and $I(X)>1$.  As explained in \cite[§2.2, Remark~2]{boctsko},
any smooth cubic projective hypersurface of dimension at least~$3$ without rational points over a $p$\nobreakdash-adic field satisfies these conditions,
thanks to a theorem of Coray and to Max Noether's theorem.
Using Theorem~\ref{mainapp} below,
it is also possible, although even more involved, to give examples of (geometrically) rational surfaces over~$\Q$
for which $\Pgen(X)=1$ and $I(X)>1$.  However such examples do not exist in dimension~$1$ (one can show that $\Pgen(X)=I(X)$ for any curve~$X$ over
any field~$k$).

The remaining of this section will be devoted to the proof of Theorem~\ref{thorder}.
We refer the reader to §\ref{sectarith} for applications of arithmetical interest.

\bigskip
\begin{proof}[ of Theorem~\ref{thorder}]%
For any dense open subset $U \subseteq X$, let $\T^0_{U/k}$ be the $k$\nobreakdash-torus whose character group is $\kbar[U]^\star/\kbar^\star$,
where~$\kbar[U]$ denotes the ring of regular functions on $U \otimes_k \kbar$.
There exist a $k$\nobreakdash-torsor $\T^1_{U/k}$ under $\T^0_{U/k}$ and a $k$\nobreakdash-morphism $U \rightarrow \T^1_{U/k}$
which satisfy the following universal property: for any $k$\nobreakdash-torus~$R^0$ and any $k$\nobreakdash-torsor~$R^1$ under~$R^0$,
any $k$\nobreakdash-morphism $U \rightarrow R^1$ factors uniquely through~$\T^1_{U/k}$ (see \cite[Lemma~2.4.4]{skobook}).

The exact sequence
of discrete $\Gamma$-modules
$$
\xymatrix{
0 \ar[r] & \kbar^\star \ar[r] & \kbar[U]^\star \ar[r] & \kbar[U]^\star / \kbar^\star \ar[r] & 0
}
$$
splits if and only if $\T^1_{U/k}(k)\neq\emptyset$ (see \cite[Lemma~2.4.4]{skobook}).
Similarly, for any $n \geq 1$,
if~$\Urond$ denotes the set of all dense open subsets of~$X$ and~$\T^n_{U/k}$ denotes the $n$\nobreakdash-fold contracted product of $\T^1_{U/k}$
as a torsor under $\T^0_{U/k}$,
then $\limproj_{U\in\Urond}(\T^n_{U/k}(k))\neq\emptyset$ if and only if there exists a family $(r_U)_{U\in\Urond}$ of
$\Gamma$-equivariant homomorphisms $r_U \colon \kbar[U]^\star \rightarrow \kbar^\star$
such that $r_U(x)=x^n$ for all $x \in \kbar^\star$ and all~$U$, and such that~$r_U$ and~$r_V$ coincide on~$\kbar[U]^\star$
for all~$U$ and all $V \subseteq U$.
Now this condition holds if and only if there exists a $\Gamma$-equivariant homomorphism $r \colon \kbar(X)^\star \rightarrow \kbar^\star$
such that $r(x)=x^n$ for all $x \in \kbar^\star$.
This in turn is equivalent to the order of~$\ob(X)$ dividing~$n$.  We are thus reduced to proving
that $\limproj_{U\in\Urond}(\T^n_{U/k}(k))\neq\emptyset$ if $n=\Pgen(X)$.

The universal property of the Albanese torsor provides us with a $k$\nobreakdash-morphism
$\Alb^n_{U/k} \rightarrow \T^n_{U/k}$ for each $U\in\Urond$,
functorially in~$U$.  Hence there is a map
$\limproj_{U\in\Urond} (\Alb^n_{U/k}(k)) \rightarrow \limproj_{U\in\Urond} (\T^n_{U/k}(k))$.
In particular we need only prove that $\limproj_{U\in\Urond}(\Alb^n_{U/k}(k))\neq\emptyset$ as soon as $\Pgen(X)$ divides~$n$, and this is achieved
by the following proposition.
\end{proof}

\bigskip
\begin{proposition}
\label{proplimprojalb}
Let~$X$ be a smooth proper geometrically integral variety over~$k$ and let~\mbox{$n\geq 1$}.
Denote by~$\Urond$ the set of all dense open subsets of~$X$.
If $\Alb^n_{U/k}(k)\neq\emptyset$ for all $U \in \Urond$, then $\limproj_{U\in\Urond}(\Alb^n_{U/k}(k))\neq\emptyset$.
\end{proposition}

\bigskip
\begin{proof}
We begin with two lemmas.

\bigskip
\begin{lemma}
\label{proplimprojalblemma1}
Let $U \subseteq V \subseteq X$ be dense open subsets, and suppose that $\NS(V \otimes_k \kbar)=0$.
Then there is an exact sequence of $k$\nobreakdash-group schemes
$$
\xymatrix{
0 \ar[r] & Q \ar[r] & \Alb^0_{U/k} \ar[r]^m & \Alb^0_{V/k} \ar[r] & 0
}
$$
where~$Q$ is a quasi-trivial $k$\nobreakdash-torus and~$m$ is induced
by the inclusion $U \subseteq V$.
\end{lemma}

\bigskip
(Recall that a $k$\nobreakdash-torus is said to be \emph{quasi-trivial} if it is isomorphic to a product of
Weil restrictions of scalars of~$\Gm$ from finite separable extensions of~$k$ down to~$k$.)

\bigskip
\begin{proof}
For any dense open subset $W \subseteq X$, the kernel of the natural map
$\Alb^0_{W/k} \rightarrow \Alb^0_{X/k}$ is $\TAlb^0_{W/k}$
(recall that $\TAlb^0_{W/k}$ denotes the largest $k$\nobreakdash-torus contained in~$\Alb^0_{W/k}$),
and the group of characters of this torus identifies with the group
$\Div^0_{X\setminus W}(X \otimes_k \kbar)$ of divisors on $X \otimes_k \kbar$ which are algebraically equivalent to~$0$
and supported on~$X \setminus W$ (see Theorem~\ref{thserreexpl}).
As a consequence, we need only prove that the map $\TAlb^0_{U/k} \rightarrow \TAlb^0_{V/k}$ induced by~$m$ is surjective
and that its kernel is a quasi-trivial $k$\nobreakdash-torus.  In other words, we need only prove that the natural map
$$
\Div^0_{X \setminus V}(X \otimes_k \kbar) \longrightarrow \Div^0_{X \setminus U}(X \otimes_k \kbar)
$$
is injective and that its cokernel is a permutation
$\Gamma$\nobreakdash-module.  Injectivity is tautological.  The hypothesis
$\NS(V \otimes_k \kbar)=0$ implies that the cokernel is naturally isomorphic
to $\Div_{V \setminus U}(V\otimes_k \kbar)$, which is indeed a permutation
$\Gamma$\nobreakdash-module.
\end{proof}

\bigskip
\begin{lemma}
\label{proplimprojalblemma2}
Let $U$, $V$ be dense open subsets of~$X$ such that $\NS((U \cup V) \otimes_k \kbar)=0$.  For any $n \geq 0$, the natural
morphism
$$
\Alb^n_{(U \cap V)/k} \longrightarrow \Alb^n_{U/k} \times_{\Alb^n_{(U \cup V)/k}} \Alb^n_{V/k}
$$
is an isomorphism.
\end{lemma}

\bigskip
\begin{proof}
We may assume $k=\kbar$ and then $n=0$.  Consider the commutative diagram of algebraic groups
$$
\xymatrix{
& 0 \ar[d] & 0 \ar[d] & 0 \ar[d] \\
0 \ar[r] & \ar[d]\TAlb^0_{(U \cap V)/k} \ar[r] & \ar[d]\TAlb^0_{U/k} \times \TAlb^0_{V/k} \ar[r] & \ar[d]\TAlb^0_{(U \cup V) /k} \\
0 \ar[r] & \ar[d]\Alb^0_{(U \cap V)/k} \ar[r] & \ar[d]\Alb^0_{U/k} \times \Alb^0_{V/k} \ar[r] & \ar[d]\Alb^0_{(U \cup V) /k} \\
0 \ar[r] & \Alb^0_{X/k} \ar[r]\ar[d] & \Alb^0_{X/k} \times \Alb^0_{X/k} \ar[r] \ar[d] & \Alb^0_{X/k} \ar[d] \rlap{\;\!\text{,}} \\
& 0 & 0 & 0
}
$$
in which the leftmost (resp.~rightmost) horizontal maps are the products (resp.~differences) of the maps induced by the various inclusions.
The statement of the lemma amounts to the exactness of the middle row.  The bottom row and the columns are exact, so it suffices
to prove that the top row is exact.  For this it suffices to prove that the corresponding sequence of character groups is exact.
By Theorem~\ref{thserreexpl}, the latter identifies with
$$
\xymatrix{
\Div^0_{F \cap G}(X) \ar[r] & \Div^0_F(X) \times \Div^0_G(X) \ar[r] & \Div^0_{F\cup G}(X) \ar[r] & 0\rlap{\text{,}}
}
$$
where $F=X \setminus U$, $G=X \setminus V$, the first map sends a divisor~$D$ to $(D,-D)$, and the second
map sends $(D,D')$ to $D+D'$.  Exactness in the middle is obvious.  Exactness on the right follows from the hypothesis that
$\NS(U \cup V)=0$.
\end{proof}

\bigskip
We now return to the proof of Proposition~\ref{proplimprojalb}.
Let $V \subseteq X$ be a dense open subset such that $\NS(V\otimes_k \kbar)=0$ (such a~$V$
exists because $\NS(X\otimes_k \kbar)$ is finitely generated).
Choose $x_V \in \Alb^n_{V/k}(k)$.  Let $F \subset V$ be an irreducible closed subset of codimension~$1$.
The natural morphism
\mbox{$m \colon \Alb^n_{(V\setminus F)/k} \rightarrow \Alb^n_{V/k}$} is a torsor under the $k$\nobreakdash-torus~$Q$ of Lemma~\ref{proplimprojalblemma1}.
Since~$Q$ is quasi-trivial, this implies that all fibres of~$m$ contain rational points, by Shapiro's lemma and Hilbert's Theorem~90.
We can therefore choose for each such~$F$ a point $x_{V \setminus F} \in \Alb^n_{(V \setminus F)/k}(k)$ whose image by~$m$ is~$x_V$.
Now let~$U$ be an arbitrary dense open subset of~$V$.  Denote by $\uplet{F_1}{F_m}$ the irreducible components of $V \setminus U$ which are of codimension~$1$
in~$V$.  Quite generally, if~$W$ is a dense open subset of~$X$ and $Z \subset W$ is a closed subset of codimension at least~$2$,
the natural map $\Alb^n_{(W\setminus Z)/k} \rightarrow \Alb^n_{W/k}$ is an isomorphism (this is a consequence of~\cite[Lemma~3.3]{milneav}).
In particular, 
the natural morphism
$\Alb^n_{U/k} \longrightarrow \Alb^n_{(V \setminus (\bigcup_{1 \leq i \leq m} F_i))/k}$ is an isomorphism,
as well as the natural morphism
$\Alb^n_{(V \setminus (\bigcap_{1 \leq i \leq m} F_i))/k} \rightarrow \Alb^n_{V/k}$ for any subset $I \subseteq \{\uplet{1}{m}\}$ of
cardinality at least~$2$.
In view of these remarks,
Lemma~\ref{proplimprojalblemma2} implies that the natural morphism
$$
\Alb^n_{U/k} \longrightarrow \Alb^n_{(V \setminus F_1)/k} \times_{\Alb^n_{V/k}} \cdots \times_{\Alb^n_{V/k}} \Alb^n_{(V \setminus F_m)/k}
$$
is an isomorphism.
Define $x_U \in \Alb^n_{U/k}(k)$ to be the inverse image of $(\uplet{x_{V\setminus F_1}}{x_{V\setminus F_m}})$.
Letting~$\Urond_V$ denote $\ensemble{U \in \Urond}{U \subseteq V}$,
we have now produced a whole family $x=(x_U)_{U \in \Urond_V} \in \prod_{U \in \Urond_V}(\Alb^n_{U/k}(k))$.
It is immediate that~$x$ is an element
of $\limproj_{U \in \Urond_V}(\Alb^n_{U/k}(k))$, hence $\limproj_{U \in \Urond}(\Alb^n_{U/k}(k))\neq\emptyset$.
\end{proof}

\bigskip
\begin{remark}
\label{laremarquerecip}
It is not known whether the statement of Proposition~\ref{proplimprojalb} remains true, even for $n=1$, if one replaces all occurrences of~$\Alb^n_{U/k}$ in it with~$\T^n_{U/k}$.
In other words, while the vanishing of~$\ob(X)$ clearly implies the existence of a $\Gamma$\nobreakdash-equivariant splitting of the exact sequence
$$
\xymatrix{
0 \ar[r] & \kbar^\star \ar[r] & \kbar[U]^\star \ar[r] & \kbar[U]^\star / \kbar^\star \ar[r] & 0
}
$$
for every dense open subset $U \subseteq X$, it is an open question whether the converse implication holds in general.
\end{remark}

\section{Arithmetical applications and fields of dimension~$\leq 1$}
\label{sectarith}

In this section we give several applications of Theorem~\ref{thorder} when~$k$
is a $p$\nobreakdash-adic field, a real closed field, or a number field.
Our main result here is Theorem~\ref{mainapp}, which answers positively a question
posed by Borovoi, Colliot-Thélène, and Skorobogatov in~\cite{boctsko} about the
elementary obstruction over number fields.  In addition, in the three cases under consideration, we address the
following deceptively simple-looking general question, which was also raised
by the aforementioned authors in~\cite{boctsko} (see \cite[§2.1]{boctsko};
there~$k$ is assumed to have characteristic~$0$):

\bigskip
\begin{question}
\label{questionext}
Let~$K/k$ be a field extension and~$X$ be a smooth proper geometrically integral $k$\nobreakdash-variety.  Let $\ob(X \otimes_k K)$
denote the class of the elementary obstruction on the $K$\nobreakdash-variety $X \otimes_k K$.  Does $\ob(X)=0$ imply $\ob(X \otimes_k K)=0$~?
\end{question}

\bigskip
We then proceed to consider the elementary obstruction and generic periods over fields of dimension~\mbox{$\leq 1$} (in the sense of Serre~\cite[II.§3]{serrecg}).
The results we obtain for fields of dimension~$\leq 1$ are in marked contrast with those for $p$\nobreakdash-adic or real closed fields
and number fields.

\subsection{Preliminaries}

We first establish two general lemmas that will be used below.

\bigskip
\begin{lemma}
\label{prelimlemmans}
Let~$k$ be a field and~$X$ be a smooth proper geometrically integral variety over~$k$.
Let $V \subseteq X$ be a dense open subset such that $\NS(V \otimes_k \kbar)=0$.
Then $P(V)=\Pgen(X)$.
\end{lemma}

\bigskip
\begin{proof}
Let $W \subseteq V$ be a dense open subset.  We need only prove that $P(V)=P(W)$.
By Lemma~\ref{proplimprojalblemma1} and Hilbert's Theorem~90, the Galois cohomology map
$
H^1(k, \Alb^0_{W/k}) \rightarrow H^1(k, \Alb^0_{V/k})
$
induced by the inclusion $W \subseteq V$ is injective.  It sends the class of $\Alb^1_{W/k}$ to the class of $\Alb^1_{V/k}$; therefore these
two classes have the same order.
\end{proof}

\bigskip
\begin{lemma}
\label{lemmaobpropagates}
Let~$k$ be a field and $f \colon X \dashrightarrow Y$ be a rational map between smooth geometrically integral varieties over~$k$.
If $\ob(X)=0$, then $\ob(Y)=0$.
\end{lemma}

\bigskip
\begin{proof}
After shrinking~$X$, we may assume~$f$ is a morphism.  We denote as usual~$\kbar$ a separable closure of~$k$
and $\Gamma=\Gal(\kbar/k)$.
Let $\xi \in Y \otimes_k \kbar$ be the image by~$f$ of the generic
point of $X \otimes_k \kbar$.
There is a natural exact sequence of discrete $\Gamma$\nobreakdash-modules
\begin{equation}
\label{applemmaex}
\xymatrix{
0 \ar[r] & \Orond_{Y \otimes_k \kbar, \;\xi}^\star \ar[r] & \kbar(Y)^\star \ar[r] & D \ar[r] & 0 \rlap{\text{,}}
}
\end{equation}
where~$D$ is the group of divisors on~$Y \otimes_k \kbar$ whose support contains~$\xi$.  Now~$D$ is a permutation Galois module,
so that the short exact sequence~(\ref{applemmaex}) splits (indeed $\Ext^1_\Gamma(D,\Orond_{Y\otimes_k\kbar, \;\xi}^\star)=0$
by Shapiro's lemma and Grothendieck's Hilbert~90 theorem).
Composing a $\Gamma$\nobreakdash-equivariant retraction of the leftmost map of~(\ref{applemmaex})
with the evaluation map $\Orond_{Y \otimes_k \kbar, \;\xi}^\star \rightarrow \kbar(\xi)^\star$ and the natural
inclusion $\kbar(\xi)^\star \subseteq \kbar(X)^\star$, we finally obtain a $\Gamma$\nobreakdash-equivariant map $\kbar(Y)^\star \rightarrow \kbar(X)^\star$
which induces the identity on~$\kbar^\star$.  The existence of a $\Gamma$\nobreakdash-equivariant retraction of the inclusion $\kbar^\star \subseteq \kbar(Y)^\star$ then
follows from the existence of a $\Gamma$\nobreakdash-equivariant retraction of the inclusion $\kbar^\star \subseteq \kbar(X)^\star$.
\end{proof}

\subsection{Local fields}

In \cite[Theorem~2.7]{boctsko}, Borovoi, Colliot-Thélène, and Skorobogatov proved that if~$X$ is a smooth proper geometrically integral variety
over a $p$\nobreakdash-adic field or a real closed field, then $\ob(X)=0$ implies $P(X)=1$.  This was a reformulation of a result of van Hamel~\cite{vanhamel}, which itself
generalised earlier work of Roquette and Lichtenbaum.  It is easy to see that the implication $\ob(X)=0 \Rightarrow P(X)=1$ fails to be
an equivalence over $p$\nobreakdash-adic or real closed fields in general --- any curve of genus~$0$ without rational points furnishes a counter-example.

The following theorem strengthens \cite[Theorem~2.7]{boctsko} and supplements it with a converse.

\bigskip
\begin{theorem}
\label{thapploc}
Let~$k$ be a $p$\nobreakdash-adic field or a real closed field and~$X$ be a smooth proper geometrically integral variety over~$k$.
Then $\ob(X)=0$ if and only if $\Pgen(X)=1$.
\end{theorem}

\bigskip
\begin{proof}
If $\Pgen(X)=1$, then $\ob(X)=0$ by Theorem~\ref{thorder}.  Now assume $\ob(X)=0$ and let $U \subseteq X$ be a dense open subset.
We need only prove that $P(U)=1$, \emph{i.e.}, that $\Alb^1_{U/k}(k)\neq\emptyset$.
The existence of a rational map
$X \dashrightarrow \Alb^1_{U/k}$ and the condition $\ob(X)=0$ imply together that $\ob(\Alb^1_{U/k})=0$ by~Lemma~\ref{lemmaobpropagates}.
From this it follows that $\Alb^1_{U/k}(k)\neq\emptyset$ by \cite[Theorem~3.2]{boctsko}.
\end{proof}

\bigskip
\begin{remark}
\label{rkvh}
If~$k$ is a $p$\nobreakdash-adic field,
it is possible to prove that the integer denoted $PsI(X)$ in~\cite{vanhamel} coincides with $\Pgen(X)$,
as van Hamel explained to the author.
Granting this, the conclusion of Theorem~\ref{thapploc} in the case of a $p$\nobreakdash-adic field
is also a consequence of \cite[Theorem~2]{vanhamel} and~\cite[Theorem~2.5]{boctsko}.
(Note that if~$k$ is a real closed field, the integers $PsI(X)$ and $\Pgen(X)$ need not be equal; see the remark after the proof of \cite[Theorem~2.6]{boctsko} for
an example where $\Pgen(X)=1$ but $PsI(X)=2$.)
\end{remark}

\bigskip
We can now give an affirmative answer to Question~\ref{questionext} in case the base field is a $p$\nobreakdash-adic field or a real closed field:

\bigskip
\begin{corollary}
\label{corquestlf}
Let~$k$ be a $p$\nobreakdash-adic field or a real closed field and~$X$ be a smooth proper geometrically integral variety over~$k$.
If $\ob(X)=0$, then $\ob(X\otimes_k K)=0$ for any field extension~$K/k$.
\end{corollary}

\bigskip
\begin{proof}
Let $V \subseteq X$ be a dense open subset such that $\NS(V \otimes_k \kbar)=0$.
Suppose $\ob(X)=0$.  By Theorem~\ref{thapploc} we then have $\Pgen(X)=1$, so that in particular $\Alb^1_{V/k}(k)\neq\emptyset$.
Now the natural map $\Alb^1_{V \otimes_k K/K} \rightarrow \Alb^1_{V/k} \otimes_k K$ is an isomorphism,
by Corollary~\ref{corformation} (note that~$k$ is perfect).  Hence $\Alb^1_{V \otimes_k K/K}(K)\neq\emptyset$,
or equivalently $P(V \otimes_k K)=1$.
Let~$\Kbar$ be an algebraic closure of~$K$.  The condition $\NS(V \otimes_k \kbar)=0$ implies $\NS(V\otimes_k \Kbar)=0$.
We can therefore apply Lemma~\ref{prelimlemmans} to the $K$\nobreakdash-variety $V \otimes_k K$ and obtain the equality $\Pgen(X \otimes_k K)=1$.
Invoking Theorem~\ref{thorder} now concludes the proof.
\end{proof}

\subsection{Number fields}

Let~$k$ be a number field and~$X$ be a smooth geometrically integral variety over~$k$.  Let~$\A_k$ denote the ring of adèles of~$k$.
In this context, obstructions to the existence of rational points
more elaborate than the elementary obstruction have been defined.  Of particular importance is the Brauer-Manin obstruction associated
with the subgroup $\Ba(X) \subseteq \Br(X)$ of locally constant algebraic classes in the cohomological Brauer group $\Br(X)=H^2_{\text{ét}}(X,\Gm)$.
This obstruction is said to vanish when \mbox{$X(\A_k)^\Ba\neq\emptyset$}, where $X(\A_k)^\Ba$ denotes the set of adelic points of~$X$
which are orthogonal to~$\Ba(X)$ with respect to the Brauer-Manin pairing $X(\A_k) \times \Br(X) \rightarrow \Q/\Z$ (\cite[§5.2]{skobook}).

Borovoi, Colliot-Thélène, and Skorobogatov proved in \cite[Theorem~2.13]{boctsko} that $\ob(X)=0$ implies $X(\A_k)^\Ba\neq \emptyset$,
provided $X(\A_k)\neq\emptyset$.  In the remark which follows their proof, they note that the converse is known to hold
by a theorem of Colliot-Thélène and Sansuc
if $\Pic(X\otimes_k\kbar)$ is a free abelian group (as is the case for instance if~$X$ is geometrically unirational);
see \cite[Prop.~3.3.2]{ctsandesc2}.
They then ask whether the converse holds in general.  The following theorem provides an affirmative answer to this question, modulo
the widely believed conjecture that Tate-Shafarevich groups of abelian varieties over number fields are finite.  It also provides
a converse to a theorem of Eriksson and Scharaschkin (\cite[Theorem~1.1]{dennis}).

\bigskip
\begin{theorem}
\label{mainapp}
Let~$k$ be a number field and~$X$ be a smooth proper geometrically integral variety over~$k$.  Assume that the Tate-Shafarevich group of the Picard variety
of~$X$ over~$k$ is finite.  Consider the following conditions:
\begin{enumerate}
\item[(i)] $X(\A_k)^\Ba\neq\emptyset$;
\item[(ii)] $\ob(X)=0$;
\item[(iii)] $\Pgen(X)=1$.
\end{enumerate}
Conditions~(ii) and~(iii) are equivalent. If~$X(\A_k)\neq\emptyset$, then~(i), (ii), and~(iii) are equivalent.
\end{theorem}

\bigskip
Note that the properness assumption on~$X$ is harmless; indeed, for any dense open subset $U \subseteq X$, the conditions
$X(\A_k)^\Ba \neq \emptyset$ and $U(\A_k)^\Ba \neq \emptyset$ are equivalent, since on the one hand
$X(\A_k)\neq\emptyset$ implies $U(\A_k)\neq\emptyset$ (by Nishimura's lemma~\cite[IV.6.2]{kollar}, the Lang-Weil estimates, and Hensel's lemma) and on the other hand $\Ba(U)=\Ba(X)$ (for a proof of this equality
see \cite[Lemme~6.1]{sansuclin} or \cite[Lemma~3.4]{dennis}).

\bigskip
\begin{proof}[ of Theorem~\ref{mainapp}]%
We have (iii)$\Rightarrow$(ii) by Theorem~\ref{thorder}, (ii)$\Rightarrow$(i) if $X(\A_k)\neq\emptyset$
by \cite[Theorem~2.13]{boctsko}, and (i)$\Rightarrow$(iii)
by \cite[Theorem~1.1]{dennis}.  This proves Theorem~\ref{mainapp} in case $X(\A_k)\neq\emptyset$.

(The finiteness of the relevant Tate-Shafarevich group is only used for establishing~(i)$\Rightarrow$(iii).
The proof of \cite[Theorem~1.1]{dennis} in our context is as follows: $X(\A_k)^\Ba\neq\emptyset$ implies $U(\A_k)^\Ba\neq\emptyset$
(where $U \subseteq X$ is an arbitrary dense open subset), which implies $\Alb^1_{U/k}(\A_k)^\Ba\neq\emptyset$ (thanks to the
existence of a morphism $U\rightarrow \Alb^1_{U/k}$), and then it follows that $\Alb^1_{U/k}(k)\neq\emptyset$ (and thus $P(U)=1$)
by a theorem of Harari and Szamuely~\cite[Theorem~1.1]{hasz}, which relies on the finiteness of the Tate-Shafarevich group of the Picard variety of~$X$.)

Let us now prove Theorem~\ref{mainapp} in general.
As Theorem~\ref{thorder} still yields (iii)$\Rightarrow$(ii), we need only
consider (ii)$\Rightarrow$(iii).  Let $U \subseteq X$ be a dense open subset such that $P(U)=\Pgen(X)$.
By Lemma~\ref{lemmaobpropagates} applied to the canonical rational map $X \dashrightarrow \Alb^1_{U/k}$,
we may assume that~$U$ is a torsor under a semi-abelian variety, after replacing~$U$ with~$\Alb^1_{U/k}$
and~$X$ with a smooth compactification of~$\Alb^1_{U/k}$.
(To ensure the existence of such a compactification, one can avoid Hironaka's theorem on resolution of singularities; see the comment after Corollary~\ref{cordoesntdepend}.)
In this case, it follows from condition~(ii) that $X(\A_k)\neq\emptyset$
by \cite[Proposition~2.12 and Theorem~3.2]{boctsko}, and
we have already proven Theorem~\ref{mainapp} under this assumption.
\end{proof}

\bigskip
The following corollary gives a conditional affirmative answer to Question~\ref{questionext} in case the base field is a number field.

\bigskip
\begin{corollary}
\label{corquestnf}
Let~$X$ be a smooth proper geometrically integral variety over a number field~$k$.  Assume that the Tate-Shafarevich group of the Picard variety of~$X$ over~$k$ is finite.
If $\ob(X)=0$, then $\ob(X\otimes_k K)=0$ for any field extension~$K/k$.
\end{corollary}

\bigskip
\begin{proof}
Corollary~\ref{corquestnf} follows from Theorem~\ref{mainapp} in the exact same way as Corollary~\ref{corquestlf} follows from Theorem~\ref{thapploc}.
\end{proof}

\subsection{Fields of dimension~$\leq 1$}

Recall that a field~$k$ is said to be a field of dimension~$\leq 1$ if the Brauer group of every finite extension of~$k$ is trivial.
For perfect fields, this is equivalent to having cohomological dimension~$\leq 1$ (\cite[II.§3.1, Proposition~6]{serrecg}).
Examples of fields of dimension~$\leq 1$ are $(C_1)$ fields (\cite[II.§3.2]{serrecg}).

Many properties of algebraic varieties over arbitrary fields have a tendency
to either become trivial or to not get any easier to understand
when they are specialised to perfect fields of dimension~$\leq 1$.
Torsors under connected linear algebraic groups thus always have rational points over such fields (as proven by Steinberg, see \cite[III.§2.3]{serrecg}),
whereas torsors under abelian varieties do not (examples are easy to construct over $\C((t))$; in particular
generic periods of varieties over perfect fields of dimension~$\leq 1$ can
be nontrivial).  In view of Theorems~\ref{thapploc} and~\ref{mainapp}, this
prompts the question whether~$\Pgen(X)=1$ might be equivalent to $\ob(X)=0$
over perfect fields of cohomological dimension~$\leq 1$.  The answer is in the negative, as shown by the following theorem (whose
proof does not make use of Theorem~\ref{thorder}).

\bigskip
\begin{theorem}
\label{thcdone}
Let~$k$ be a field of dimension~$\leq 1$.  Then $\ob(X)=0$ for any geometrically integral variety~$X$ over~$k$.
\end{theorem}

\bigskip
We note that if~$k$ is a field of dimension~$\leq 1$, the exact sequence of Galois modules
$$
\xymatrix{
0 \ar[r] & \kbar^\star \ar[r] & \kbar[U]^\star \ar[r] & \kbar[U]^\star/\kbar^\star \ar[r] & 0
}
$$
is clearly split
for any dense open subset $U \subseteq X$ (indeed the torsor~$\T^1_{U/k}$
which appears in the proof of Theorem~\ref{thorder} has a rational point since it is a torsor under a torus and~$k$
has dimension~$\leq 1$ \cite[p.~170]{corpslocaux}).  If this implied $\ob(X)=0$, Theorem~\ref{thcdone} would follow immediately.
See Remark~\ref{laremarquerecip}.

\bigskip
To prove Theorem~\ref{thcdone} we start with a variant of a well-known lemma in homological algebra originally due to Matlis (\cite[Lemma~1]{matlis}).

\bigskip
\begin{lemma}
\label{lemmamatlis}
Let~$G$ be a profinite group and~$M$ be a discrete $G$\nobreakdash-module.
Fix $n \geq 1$ and assume $\Ext^n_G(\Z[G/H]/I,M)=0$ for any normal open subgroup $H \subseteq G$
and any left ideal $I \subseteq \Z[G/H]$.  Then $\Ext^n_G(B,M)=0$ for any discrete $G$\nobreakdash-module~$B$.
\end{lemma}

\bigskip
\begin{proof}
We prove the lemma by induction on~$n$.  If the assertion is true for a given $n \geq 1$, it is seen to hold
for $n+1$ by embedding~$M$ into an injective discrete $G$\nobreakdash-module~$M'$ and noting that
$\Ext^{n+1}_G(B,M)=\Ext^{n}_G(B,M'/M)$.  Thus we can assume $n=1$, and the conclusion of the lemma is then equivalent
to~$M$ being injective.  Let $E \subseteq F$ be an inclusion of
discrete $G$\nobreakdash-modules and $f \colon E \rightarrow M$ be a $G$\nobreakdash-equivariant morphism.
By Zorn's lemma, there exists a sub\nobreakdash-$G$\nobreakdash-module $E' \subseteq F$ containing~$E$ and
a $G$\nobreakdash-equivariant morphism $f' \colon E' \rightarrow M$ which extends~$f$ and which cannot be
extended to a $G$\nobreakdash-equivariant morphism on any larger sub\nobreakdash-$G$\nobreak-module of~$F$.
It only remains to be shown that $F=E'$.
Let $x \in F$ and let $H \subseteq G$ be a normal open subgroup contained in the stabiliser of~$x$.
Let $E'' \subseteq F$ denote the sub\nobreakdash-$G$\nobreakdash-module of~$F$ generated by~$E'$ and~$x$.
Then $E''/E'$ is of the shape $\Z[G/H]/I$ for some left ideal~$I$, and as a consequence we have $\Ext^1_G(E''/E',M)=0$.
In particular~$f$ can be extended to a $G$\nobreakdash-equivariant morphism $E'' \rightarrow M$, so that $E''=E'$
and therefore $x \in E'$; hence finally $F=E'$.
\end{proof}

\bigskip
\begin{proof}[ of Theorem~\ref{thcdone}]%
Let~$\kbar$ be a separable closure of~$k$ and $\Gamma=\Gal(\kbar/k)$.
According to \cite[Proposition~2.2.4]{ctsandesc2}, it suffices to establish that $\Ext^2_\Gamma(\Pic(X \otimes_k \kbar),\kbar^\star)=0$.
We shall actually prove that $\Ext^2_\Gamma(B,\kbar^\star)=0$ for any discrete $\Gamma$\nobreakdash-module~$B$.
By Lemma~\ref{lemmamatlis}, we may assume that $B=\Z[\Gamma/H]/I$, where $H \subseteq \Gamma$ is a normal open subgroup
and $I \subseteq \Z[\Gamma/H]$ is a left ideal.
It follows from the exact sequence
$$
\xymatrix{
0 \ar[r] & I \ar[r] & \Z[\Gamma/H] \ar[r] & B \ar[r] & 0
}
$$
that in order for $\Ext^2_\Gamma(B,\kbar^\star)$ to be trivial, it suffices
that $\Ext^2_\Gamma(\Z[\Gamma/H],\kbar^\star)$ and $\Ext^1_\Gamma(I,\kbar^\star)$ be trivial.
The first of these groups is equal to the Brauer group of a finite extension of~$k$ (quite generally
one has $\Ext^n_\Gamma(\Z[\Gamma/H],-)=H^n(H,-)$ for $n=0$ hence for any~$n$), therefore it is trivial.
By \cite[Ch.~I, §0, Theorem~0.3 and Example~0.8]{milneadt}, which we can apply because~$I$ is finitely generated as a $\Z$\nobreakdash-module,
the second of these groups fits into an exact sequence
$$
\xymatrix{
0 \ar[r] & H^1(k, \Hom_\Z(I,\kbar^\star)) \ar[r] & \Ext^1_\Gamma(I,\kbar^\star) \ar[r] & \Ext^1_\Z(I,\kbar^\star) \rlap{\text{,}}
}
$$
where the subscript~$\Z$ refers to the category of abelian groups.  Now~$I$ is a free $\Z$\nobreakdash-module;
hence on the one hand, the group $\Ext^1_\Z(I,\kbar^\star)$ is trivial, and on the other hand,
the group $H^1(k, \Hom_\Z(I,\kbar^\star))$ is equal to $H^1(k,T)$ where~$T$ is a $k$\nobreakdash-torus,
and is therefore also trivial as~$k$ is a field of dimension~$\leq 1$ (see~\cite[p.~170]{corpslocaux}).
\end{proof}

\bigskip
In the statement of Theorem~\ref{thcdone}, the hypothesis that~$k$ be a field of dimension~$\leq 1$ cannot be replaced with the weaker hypothesis that~$k$ have cohomological
dimension~$\leq 1$.  In other words, there are geometrically integral varieties defined over (imperfect) fields of cohomological dimension~$1$
for which the elementary obstruction does not vanish.
Even better, the statement of Theorem~\ref{thcdone} turns out to be optimal:

\bigskip
\begin{proposition}
Let~$k$ be a field.  Suppose that $\ob(X)=0$ for any geometrically integral variety~$X$ over~$k$.  Then~$k$ is a field of dimension~$\leq 1$.
\end{proposition}

\bigskip
\begin{proof}
It suffices to prove that $\Br(K)=0$ for all finite separable extensions~$K/k$ (\cite[II.§3.1, Proposition~5]{serrecg}).
Let~$X$ be a Severi-Brauer variety over~$K$.  It is a toric $K$\nobreakdash-variety (see for instance \cite[Example~8.5.1]{merkpan});
as a consequence, the Weil restriction of scalars $R_{K/k}X$ is a toric $k$\nobreakdash-variety.
Now we have $\ob(R_{K/k}X)=0$ by hypothesis, so that $(R_{K/k}X)(k)\neq\emptyset$ (see \cite[Lemma~2.1~(iv)]{boctsko}, whose proof works for
arbitrary fields) and therefore $X(K)\neq\emptyset$.
\end{proof}

\bigskip
Let us now consider Question~\ref{questionext} under the assumption that~$k$ is a field of dimension~$\leq 1$.
An affirmative answer would imply,
together with Theorem~\ref{thcdone}, that $\ob(X \otimes_k K)=0$
for any geometrically integral $k$\nobreakdash-variety~$X$ and any field extension~$K/k$.
We only establish the following weaker result (whose proof does make use of Theorem~\ref{thorder}).

\bigskip
\begin{proposition}
\label{propcdoneext}
Let~$k$ be a perfect field of dimension~$\leq 1$ and~$X$ be a smooth proper geometrically integral variety over~$k$.
If $P(X)=1$, then $\ob(X \otimes_k K)=0$ for any field extension $K/k$.
\end{proposition}

\bigskip
\begin{proof}
Let $V \subseteq X$ be a dense open subset and $n \geq 0$.  The natural morphism $\Alb^n_{V/k} \rightarrow \Alb^n_{X/k}$ is a torsor
under $\TAlb^0_{V/k}$.  In particular its fibre over any $k$\nobreakdash-point of $\Alb^n_{X/k}$
is a $k$\nobreakdash-torsor under a torus, hence contains a rational point by \cite[p.~170]{corpslocaux}.  This proves that $\Pgen(X)=1$.
Now the rest of the argument is the same as in the proof of Corollary~\ref{corquestlf}.
\end{proof}

\section{Albanese torsors, $1$-motives, and the elementary obstruction}
\label{secquat}

Let~$X$ be a smooth proper geometrically
integral variety over a field~$k$.  As is well known, the elementary obstruction to the existence of a $0$\nobreakdash-cycle of degree~$1$ on~$X$
vanishes if and only if the $2$\nobreakdash-extension~$e(X)$ of $\Pic(X \otimes_k \kbar)$ by $\kbar^\star$ defined by the exact
sequence of $\Gamma$\nobreakdash-modules
$$
\xymatrix{
0 \ar[r] & \kbar^\star \ar[r] & \kbar(X)^\star \ar[r] & \Div(X\otimes_k \kbar) \ar[r] & \Pic(X \otimes_k \kbar) \ar[r] & 0
}
$$
gives rise to the trivial class in $\Ext^2_\Gamma(\Pic(X\otimes_k
\kbar),\kbar^\star)$, where~$\kbar$ denotes a separable closure of~$k$ and $\Gamma=\Gal(\kbar/k)$
(\cite[Proposition~2.2.4]{ctsandesc2}).  More generally, the order
of~$\ob(X)$ is equal to the order of the class of~$e(X)$ in $\Ext^2_\Gamma(\Pic(X\otimes_k \kbar),\kbar^\star)$.
Theorem~\ref{thorder} can thus be reformulated as stating that the order of the class of~$e(X)$ divides the generic period~$\Pgen(X)$;
this divisibility is not an equality in general, according to Theorem~\ref{thcdone}.

Defining an analogous $2$\nobreakdash-extension~$E(X)$ of the relative Picard
functor~$\Pic_{X/k}$ by the multiplicative group~$\Gm$ in the category of étale sheaves on~$\Sm/k$ is a straightforward task
(see §\ref{subsecmotivic} for details; $\Sm/k$ denotes the category of smooth $k$\nobreakdash-schemes of finite type).
One of our main results below is the following: the order of the class of~$E(X)$ in $\Ext^2(\Pic_{X/k},\Gm)$
is actually \emph{equal} to~$\Pgen(X)$ (see Corollary~\ref{quatcor}).
This can be seen as an optimal strengthening of Theorem~\ref{thorder},
since the order of the class of~$e(X)$, which is equal to the order of~$\ob(X)$, clearly divides the order of the class of~$E(X)$.
Our goal, however, is really to uncover a deeper connection between Albanese
torsors and the elementary obstruction than a mere comparison between the orders of two classes.
Let~$\EM(U)$ denote the $2$\nobreakdash-extension of the Picard $1$\nobreakdash-motive of~$U$ by~$\Gm$ obtained by pullback from~$E(X)$.
The main theorem of this section (Theorem~\ref{themaintheorem}) states that by applying a certain canonical construction to the $2$\nobreakdash-extension~$\EM(U)$,
one can recover the torsor~$\Alb^1_{U/k}$ as a sheaf on~$\Sm/k$ (and hence as a variety, since $\Alb^1_{U/k} \in \Sm/k$).
As mentioned in the introduction, it will follow that rational points of~$\Alb^1_{U/k}$ correspond naturally to ``Yoneda trivialisations'' of the $2$\nobreakdash-extension~$\EM(U)$.
In the case where $U=X$,
a Yoneda trivialisation of~$\EM(X)$ is simply a zigzag of morphisms of $2$\nobreakdash-extensions of~$\Pic^0_{X/k}$
by~$\Gm$, which starts with the pullback of~$\E(X)$ by the inclusion $\Pic^0_{X/k} \rightarrow \Pic_{X/k}$, and ends with the trivial $2$\nobreakdash-extension
$$
\xymatrix{
0 \ar[r] & \Gm \ar[r]^{\Id} & \Gm \ar[r]^0 & \Pic^0_{X/k} \ar[r]^{\Id} & \Pic^0_{X/k} \ar[r] & 0\rlap{\text{.}}
}
$$
When~$X$ is a torsor under an abelian variety, Theorem~\ref{themaintheorem} identifies~$X(k)$ with the set of such zigzags up to a certain homotopy relation (and similarly
for any finite separable extension $\ell/k$).

The remainder of the paper is organised as follows. In §\ref{subsecmotivic} we define the $2$\nobreakdash-extensions~$\E(X)$ and~$\EM(U)$,
and establish some of their properties.  Then in~§\ref{subsecmaincons} we describe the general construction which allows one to
recover~$\Alb^1_{U/k}$ from~$\EM(U)$, state a key theorem (Theorem~\ref{mainconsthm}), and deduce Theorem~\ref{themaintheorem} from it.
Section~§\ref{subsecexplicit} is devoted to an independent result about Poincaré sheaves on abelian varieties.  The proof of
Theorem~\ref{mainconsthm} is finally given in~§\ref{subsecproof}. It rests on the contents of~§\ref{subsecexplicit}.

Under the additional hypothesis that the field~$k$ has characteristic~$0$, another approach to some of the results contained in this section
was considered earlier by van Hamel in~\cite{vanhamel}, based on the well-known explicit description of the Weil pairing.  Specifically,
the first assertion of Corollary~\ref{quatcor} below, in the case where $U=X$ and~$k$ has characteristic~$0$, follows from \cite[Theorem~3.6]{vanhamel}
(as well as from \cite[Proposition~2.1]{skoronote}).  See also Remark~\ref{rkvh} and \cite[Remark~1.7]{vanhamellms}.

\subsection{Motivic $2$-extensions}
\label{subsecmotivic}

We fix once and for all a field~$k$ and a smooth proper geometrically integral variety~$X$ over~$k$.  (The variety~$X$ must be normal for Lemma~\ref{lemmabigecomplex} below to hold; but it is mostly
for convenience that we moreover assume it to be smooth.)

If~$S$ is a scheme, let $\Sm/S$ denote the category of smooth $S$\nobreakdash-schemes of finite type
and  $\Sh(\Sm/S)$ the category of sheaves in abelian groups on $\Sm/S$ for the étale topology.
For $S \in \Sm/k$ and $Y \in \Sm/S$, we denote by~$\Div(Y/S)$ the subgroup of~$\Div(Y)$ consisting of those divisors whose support does not
contain any irreducible component of any fibre of $Y \rightarrow S$, we denote by~$\R(Y/S)^\star$ the group of invertible rational functions on~$Y$
whose divisor belongs to~$\Div(Y/S)$, and we let
$\Pic(Y/S)$ denote the cokernel of the natural map $\Pic(S) \rightarrow \Pic(Y)$.  We let $\Pic^0(Y/S)$ be the subgroup of $\Pic(Y/S)$ consisting
of those classes whose restriction to every geometric fibre of $Y \rightarrow S$ is algebraically equivalent to~$0$,
and $\Div^0(Y/S)$ be the inverse image of $\Pic^0(Y/S)$ by the natural map $\Div(Y/S)\rightarrow\Pic(Y/S)$.
If~$L \subseteq Y$ is a subscheme, we denote by $\Div(Y/S,L)$ the subgroup of~$\Div(Y/S)$ consisting of those divisors
whose support is contained in~$L$.
For $T \in \Sm/S$, the groups $\R(Y \times_S T/T)^\star$, $\Div(Y \times_S T/T)$, $\Pic(Y \times_S T/T)$,
$\Div^0(Y \times_S T/T)$, $\Pic^0(Y \times_S T/T)$, and $\Div(Y \times_S T/T, L \times_S T)$
depend functorially on~$T$; as a consequence they define presheaves on $\Sm/S$.
Let $\Rstar_{Y/S}$, $\Div_{Y/S}$, $\Pic_{Y/S}$, $\Div^0_{Y/S}$, $\Pic^0_{Y/S}$, and $\Div_{Y/S,L}$ denote the associated étale sheaves,
and let $\Rstar_{Y/S,L}$ (resp.~$\Div^0_{Y/S,L}$) be the inverse
image of $\Div_{Y/S,L}$ by the natural map $\Rstar_{Y/S} \rightarrow \Div_{Y/S}$ (resp.~$\Div^0_{Y/S} \rightarrow \Div_{Y/S}$).
If~$Y$ is proper over~$S$, then $\Pic_{Y/S}$ coincides with the relative Picard functor (see \cite[p.~201]{blr}).

\bigskip
\begin{remark}
Under our assumptions on~$S$ and~$Y$, the groups $\R(Y/S)^\star$ and $\Div(Y/S)$ respectively coincide with the groups $\mathrm{M}(Y/S)^\star$ and $\Div(Y/S)$
defined by Grothendieck in \cite[20.6.1 and 21.15.2]{ega44}.
Note that for $T \in \Sm/S$, the group of invertible rational functions on $Y \times_S T$ does not depend
functorially on~$T$, even if $S=Y=\Spec(k)$; hence the need to consider relative rational functions and relative divisors.
\end{remark}

\bigskip
\begin{lemma}
\label{lemmabigecomplex}
The natural complex of étale sheaves in abelian groups on~$\Sm/k$
\begin{equation}
\label{bigecomplex}
\xymatrix{
0 \ar[r] & \Gm \ar[r] & \Rstar_{X/k} \ar[r] & \Div_{X/k} \ar[r] & \Pic_{X/k} \ar[r] & 0
}
\end{equation}
is exact.
\end{lemma}

\bigskip
\begin{proof}
Only exactness at $\Pic_{X/k}$ requires a proof.
It suffices to establish that for any connected $T \in \Sm/k$, any $t \in T$, and any $D \in \Div(X \times_k T)$, there is an open subset~$U$
of~$T$ containing~$t$ such that the restriction of~$D$ to $X \times_k U$ is linearly equivalent, on $X \times_k U$, to an element of $\Div(X \times_k U / U)$.
Let~$\xi$ denote the generic point of the fibre of $X \times_k T \rightarrow T$ above~$t$,
and let~$\Orond$ denote the local ring of~$X \times_k T$ at~$\xi$.
The Picard group of any local ring vanishes, hence $\Pic(\Orond)=0$, and in particular the image of~$D$ in $\Pic(\Orond)$ is trivial,
which means that after replacing~$D$ with $D+P$ for some principal divisor~$P$ on $X \times_k T$, we may assume that the support of~$D$
does not contain~$\xi$.  The set~$U$ of points $u \in T$ such that the support of~$D$ does not contain the fibre of $X \times_k T \rightarrow T$
above~$u$ is an open subset of~$T$, by Chevalley's semi-continuity theorem~\cite[13.1.3]{ega43}. By assumption it contains~$t$, and the restriction of~$D$
to~$X \times_k U$ is indeed an element of $\Div(X\times_k U/U)$.
\end{proof}

\bigskip
The complex~(\ref{bigecomplex}) thus defines a $2$\nobreakdash-extension of the relative Picard functor~$\Pic_{X/k}$
by the multiplicative group~$\Gm$ in the category $\Sh(\Sm/k)$.  We let~$E(X)$ denote this $2$\nobreakdash-extension.

Let $U \subseteq X$ be a dense open subset.
With~$U$ are naturally associated two $1$\nobreakdash-motives:
the \emph{Albanese $1$\nobreakdash-motive} $M_1(U)$ and the \emph{Picard $1$\nobreakdash-motive} $M^1(U)$.
By definition $M_1(U)$ is the semi-abelian variety $\Alb^0_{U/k}$ and $M^1(U)$ is the Cartier
dual of $M_1(U)$ (see \cite[§2.1]{ramachandran}; for the definition of Cartier duality between $1$\nobreakdash-motives, see \cite[§2.4.1]{raynaudmonodr}).
As a two-term complex of sheaves on~$\Sm/k$, the Picard $1$\nobreakdash-motive
$M^1(U)$ is canonically isomorphic to $[\Div^0_{X/k,F} \rightarrow \Pic^0_{X/k}]$, where $F=X \setminus U$ (see Theorem~\ref{thserreexpl}; note that
the sheaf on $\Sm/k$ defined by the Picard variety $\Pic^0_{X/k,\mathrm{red}}$ is equal to $\Pic^0_{X/k}$).
We take the convention that as complexes of sheaves, $1$\nobreakdash-motives are concentrated in degrees~$-1$ and~$0$ and written either vertically or between brackets.
The natural diagram
$$
\xymatrix{
&&& \Div^0_{X/k,F} \ar[d] \ar@{=}[r] & \Div^0_{X/k,F} \ar[d] \\
0 \ar[r] & \Gm \ar[r] & \Rstar_{X/k} \ar[r] & \Div^0_{X/k} \ar[r] & \Pic^0_{X/k} \ar[r] & 0
}
$$
defines a $2$\nobreakdash-extension of~$M^1(U)$ by~$\Gm$ in the category~$C(\Sh(\Sm/k))$ of complexes of étale sheaves on $\Sm/k$,
where~$\Gm$ and~$\Rstar_{X/k}$ are regarded as complexes concentrated in degree~$0$.
We let~$\EM(U)$ denote this $2$\nobreakdash-extension, considered as an object of~$C(C(\Sh(\Sm/k)))$.
The notation~$\EM(U)$ is justified by Corollary~\ref{cordoesntdepend} below.

Recall that if $S$ is a $k$\nobreakdash-scheme, an $S$\nobreakdash-morphism $f \colon U \times_k S \rightarrow U' \times_k S$ is said to be universally dominant if its image
contains a dense open subset of every fibre of $U' \times_k S \rightarrow S$.
For fixed $S \in \Sm/k$, we shall now prove (in Proposition~\ref{propemufunct} below) that as objects of the category $C(C(\Sh(\Sm/S)))$,
the $2$\nobreakdash-extensions $\EM(U)|_S$ are functorial in~$U$
with respect to universally dominant $S$\nobreakdash-morphisms,
where~$-|_S$ denotes restriction to~$\Sm/S$.

Let~$X'$ be another smooth proper geometrically integral variety over~$k$ and $U' \subseteq X'$ be a dense open subset.
Let $S \in \Sm/k$ and let $f \colon U \times_k S \rightarrow U' \times_k S$ be a universally dominant $S$\nobreakdash-morphism.
Clearly~$f$ induces a morphism of sheaves $\Rstar_{X' \times_k S/S} \rightarrow \Rstar_{X \times_k S/S}$.
Let $W \subseteq X \times_k S$ denote the domain of definition of the rational map $X \times_k S \dashrightarrow X' \times_k S$
induced by~$f$.  For any $T \in \Sm/S$,
all points of codimension~$1$ of $X\times_k T$ are contained in $W \times_S T$, since $X' \times_k S \rightarrow S$ is proper
and $X \times_k S$ is normal.
Hence \mbox{$\Div(X \times_k T/T)=\Div(W \times_S T/T)$} for all $T \in \Sm/S$, and consequently $\Div_{X \times_k S/S} = \Div_{W/S}$.
Composing this canonical isomorphism with the morphism $\Div_{X'\times_k S/S} \rightarrow \Div_{W/S}$ induced by
$f \colon W \rightarrow X'\times_k S$ yields a morphism $\Div_{X' \times_k S/S} \rightarrow \Div_{X \times_k S/S}$.
We will refer to it simply as the morphism induced by~$f$.
The two morphisms $\Rstar_{X' \times_k S/S} \rightarrow \Rstar_{X \times_k S/S}$ and $\Div_{X' \times_k S/S} \rightarrow
\Div_{X\times_k S/S}$ induced by~$f$ give rise to
a natural morphism $M^1(f) \colon M^1(U')|_S \rightarrow M^1(U)|_S$ in $C(\Sh(\Sm/S))$
and more generally to a natural morphism $\EM(U')|_S \rightarrow \EM(U)|_S$ in~$C(C(\Sh(\Sm/S)))$.

\bigskip
\begin{proposition}
\label{propemufunct}
Let~$X'$, $X''$ be smooth proper geometrically integral varieties over~$k$ and $U' \subseteq X'$, $U'' \subseteq X''$ be dense open subsets.
Let $S \in \Sm/k$.
Let $f \colon U \times_k S \rightarrow U' \times_k S$ and $g \colon U' \times_k S \rightarrow U''\times_k S$ be universally dominant $S$\nobreakdash-morphisms.
Then the morphism $\EM(U'')|_S \rightarrow \EM(U)|_S$ induced by~$g \circ f$ is equal to the composition of the morphisms
$\EM(U'')|_S \rightarrow \EM(U')|_S$ and $\EM(U')|_S \rightarrow \EM(U)|_S$ induced by~$g$ and~$f$.
\end{proposition}

\bigskip
\begin{proof}
It suffices to prove that the morphism of sheaves $\Div^0_{X'' \times_k S/S} \rightarrow \Div^0_{X \times_k S/S}$ induced by~$g \circ f$ is equal to the composition of the morphisms $\Div^0_{X''
\times_k S/S} \rightarrow \Div^0_{X'\times_k S/S}$ and $\Div^0_{X' \times_k S/S} \rightarrow \Div^0_{X \times_k S/S}$ induced by~$g$ and~$f$.  Lest the
reader think that such an assertion must be trivial, we point out that the corresponding statement with~$\Div$ instead of~$\Div^0$ is clearly false, even when $S=\Spec(k)$
(for instance choose a point $P \in \P^2(k)$, let $U=U'=U''=\P^2_k \setminus \{P\}$, $f=g=\Id$, $X'=\P^2_k$, and let~$X''$ and~$X$ be the blowing-up of~$P$ in~$\P^2_k$).

It is even enough to prove that for any $T \in \Sm/S$, the map $\Div^0(X'' \times_k T/T) \rightarrow \Div^0(X \times_k T/T)$ induced by~$g\circ f$
is equal to the composition of the maps induced by~$g$ and~$f$.
Furthermore we may assume that $T=S=\Spec(k)$, by extending scalars from~$k$ to the function fields of the connected components of~$T$.
Indeed, if~$T$ is connected,
there is for $Y \in \Sm/T$ a canonical inclusion $\Div(Y \times_k T/T)\subseteq\Div(Y \times_k \eta/\eta)$ which
is functorial in~$Y$, where~$\eta$ denotes the generic point of~$T$.
Finally, we may assume that~$k$ is algebraically closed.

We are now reduced to showing that if $f \colon X \dashrightarrow X'$ denotes a dominant
rational map between smooth proper varieties over an algebraically closed
field~$k$, and $W \subseteq X$ is its domain of definition, then the formation
of the morphism $\Div^0(X') \rightarrow \Div^0(W)=\Div^0(X)$ induced by~$f$ is compatible
with composition.  Let $U \subseteq W$ and $U' \subseteq X'$ be dense open subsets such that $f(U)\subseteq U'$.
Denote $F=X \setminus U$ and $F'=X'\setminus U'$.
The morphism $\Alb^0_{U/k} \rightarrow \Alb^0_{U'/k}$ determined by $f \colon U \rightarrow U'$
induces a morphism from the character group of~$\TAlb^0_{U'/k}$ to the character group of~$\TAlb^0_{U/k}$.
These groups are canonically isomorphic to $\Div^0_{F'}(X')$ and $\Div^0_F(X)$ respectively
(Theorem~\ref{thserreexpl}).
It is straightforward to check that the morphism $\Div^0_{F'}(X') \rightarrow \Div^0_F(X)$ thus obtained
coincides with the restriction of the morphism $\Div^0(X') \rightarrow \Div^0(X)$ induced by~$f$.
Now the formation of the morphism $\TAlb^0_{U/k} \rightarrow \TAlb^0_{U'/k}$ determined by~$f$ is clearly compatible with composition.
This remark concludes the proof, thanks to the formulas $\Div^0(X)=\limind \Div^0_F(X)$
and $\Div^0(X')=\limind \Div^0_{F'}(X')$, where~$F$ (resp.~$F'$) ranges over all closed subsets of~$X$ (resp.~$X'$) of codimension~$1$.
\end{proof}

\bigskip
\begin{corollary}
\label{coremufunct}
In the situation of Proposition~\ref{propemufunct}, if~$f$ is an isomorphism, then the morphism $\EM(U')|_S \rightarrow \EM(U)|_S$ induced by~$f$
is an isomorphism.
\end{corollary}

\bigskip
On taking $S=\Spec(k)$, $U'=U$, and $f=\Id$ in Corollary~\ref{coremufunct}, we obtain:

\bigskip
\begin{corollary}
\label{cordoesntdepend}
The $2$\nobreakdash-extension~$\EM(U)$ depends on~$U$ but does not depend on the smooth compactification $U \subseteq X$
(up to a canonical isomorphism in $C(C(\Sh(\Sm/k)))$).
\end{corollary}

\bigskip
In view of Corollary~\ref{cordoesntdepend}, we may omit to specify a smooth compactification of~$U$ when talking about~$\EM(U)$, provided
that at least one such compactification exists.  For future use, we note that torsors
under semi-abelian varieties are known to admit smooth compactifications
over arbitrary fields.  This follows from the existence of smooth equivariant compactifications of tori (for which see~\cite{cthaskotore}).

The $2$\nobreakdash-extensions $\EM(-)$ are also functorial in a weaker sense with respect to morphisms which are not necessarily dominant:

\bigskip
\begin{proposition}
\label{propyondeq}
Let $f \colon V \rightarrow W$ be a morphism of smooth geometrically integral $k$\nobreakdash-varieties which admit smooth compactifications.
The $2$\nobreakdash-extension~$\EM(W)$, regarded as a $2$\nobreakdash-extension in the abelian category $C(\Sh(\Sm/k))$, is
canonically Yoneda equivalent to the pullback of~$\EM(V)$ by the morphism $M^1(W) \rightarrow M^1(V)$ induced by~$f$.
\end{proposition}

\bigskip
\begin{proof}
Let $V \subseteq Y$ and $W \subseteq Z$ denote smooth compactifications of~$V$ and~$W$.
Let $D \subseteq \Div^0_{Z/k}$ be the smallest subsheaf in groups containing
$\Div^0_{Z/k,Z\setminus\{f(v)\}}$ for all closed points $v \in V$.
Let $\Rstar \subseteq \Rstar_{Z/k}$ denote the inverse image of~$D$ by the
natural map $\Rstar_{Z/k}\rightarrow \Div^0_{Z/k}$, and let~$H$
be the complex of sheaves $[\Div^0_{Z/k,Z\setminus W} \rightarrow D]$ with~$D$ in degree~$0$.
By choosing a closed point of~$V$ and using the fact that the Picard group of a local ring vanishes,
one sees that the natural complex
\begin{equation}
\label{seyondeq}
\xymatrix{
0 \ar[r] & \Gm \ar[r] & \Rstar \ar[r] & H \ar[r] & M^1(W) \ar[r] & 0
}
\end{equation}
is an exact sequence in the abelian category $C(\Sh(\Sm/k))$.
The rational map $Y \dashrightarrow Z$ induced by~$f$ is defined on an open subset $Y' \subseteq Y$ whose complement has
codimension~\mbox{$\geq 2$} in~$Y$.  Therefore we can define natural pullback morphisms
$D \rightarrow \Div^0_{Y'/k}=\Div^0_{Y/k}$
and $\Rstar \rightarrow \Rstar_{Y/k}$, and hence a natural
morphism of exact sequences from~(\ref{seyondeq}) to~$\EM(V)$ whose leftmost (resp.~rightmost) component is the
identity of~$\Gm$ (resp.~is the
morphism $M^1(f) \colon M^1(W) \rightarrow M^1(V)$ induced by~$f$).  As a consequence, the $2$\nobreakdash-extension~(\ref{seyondeq})
is canonically Yoneda equivalent to the pullback of~$\EM(V)$ by~$M^1(f)$.
Now there is a natural morphism of exact sequences from~(\ref{seyondeq})
to $\EM(W)$ whose leftmost (resp.~rightmost) component is the identity of~$\Gm$ (resp.~of~$M^1(W)$),
so the $2$\nobreakdash-extensions $\EM(W)$ and~(\ref{seyondeq}) are also canonically Yoneda equivalent.
\end{proof}

\bigskip
Let us now turn our attention to the classes of~$\E(X)$ and of~$\EM(U)$ in the corresponding groups of extensions.
For any $S \in \Sm/k$, any $2$\nobreakdash-extension
\begin{equation}
\label{eqobjcatpre}
\xymatrix{
0 \ar[r] & \GmS \ar[r] & E \ar[r]^f & F \ar[r] & M^1(U)|_S \ar[r] & 0
}
\end{equation}
in~$C(\Sh(\Sm/S))$ induces a canonical quasi-isomorphism $C(\GmS[1] \rightarrow C(f)) \rightarrow M^1(U)|_S$,
where~$C(-)$ denotes the mapping cone of a morphism, and therefore a distinguished triangle
\begin{equation}
\label{eqdisttri}
\xymatrix{
\GmS[1] \ar[r] & C(f) \ar[r] & M^1(U)|_S \ar[r] & \GmS[2]
}
\end{equation}
in the bounded derived category $D^b(\Sh(\Sm/S))$ if~$E$ and~$F$ are bounded.
We shall speak of the class in the hyperext group $$\Ext^2_{\Sh(\Sm/S)}(M^1(U)|_S,\GmS)=\Hom_{D^b(\Sh(\Sm/S))}(M^1(U)|_S,\GmS[2])$$
determined by the $2$\nobreakdash-extension~(\ref{eqobjcatpre}) to refer to the morphism $M^1(U)|_S \rightarrow \GmS[2]$ which appears
in the exact triangle~(\ref{eqdisttri}) associated with that $2$\nobreakdash-extension.

\bigskip
\begin{proposition}
\label{propcompareemuex}
The order of the class of~$\EM(U)$ in the hyperext group $\Ext^2_{\Sh(\Sm/k)}(M^1(U),\Gm)$
divides the order of the class of~$\E(X)$ in $\Ext^2_{\Sh(\Sm/k)}(\Pic_{X/k},\Gm)$.
Let~$\kbar$ denote a separable closure of~$k$.
If $\NS(U \otimes_k \kbar)=0$,
then these orders are equal.
\end{proposition}

\bigskip
\begin{proof}
Let $Y \in C(\Sh(\Sm/k))$ denote the complex $[\Div_{X/k,F} \rightarrow \Pic_{X/k}]$, where~$\Pic_{X/k}$ is placed in degree~$0$,
and let~$B$ denote the $2$\nobreakdash-extension of~$Y$ by~$\Gm$ defined by the diagram
$$
\xymatrix{
&&& \Div_{X/k,F} \ar[d] \ar@{=}[r] & \Div_{X/k,F} \ar[d] \\
0 \ar[r] & \Gm \ar[r] & \Rstar_{X/k} \ar[r] & \Div_{X/k} \ar[r] & \Pic_{X/k} \ar[r] & 0\rlap{\text{.}}
}
$$
The natural exact sequence
$$
\xymatrix{
0 \ar[r] & \Pic_{X/k} \ar[r] & Y \ar[r] & \Div_{X/k,F}[1] \ar[r] & 0
}
$$
in~$C(\Sh(\Sm/k))$ induces a distinguished triangle
$$
\xymatrix{
\Pic_{X/k} \ar[r] & Y \ar[r] & \Div_{X/k,F}[1] \ar[r] & \Pic_{X/k}[1]
}
$$
in~$D^b(\Sh(\Sm/k))$.  Let us apply the functor $\Hom(-,\Gm[2])$ to it.  We obtain an exact sequence
$$
\xymatrix{
\Ext^1_{\Sh(\Sm/k)}(\Div_{X/k,F},\Gm) \ar[r] & \Ext^2_{\Sh(\Sm/k)}(Y,\Gm) \ar[r] & \Ext^2_{\Sh(\Sm/k)}(\Pic_{X/k},\Gm) \rlap{\text{.}}
}
$$
Its first term vanishes, by Shapiro's lemma and Hilbert's Theorem~90.
Consequently the rightmost map of this exact sequence is injective.  This map sends the class of~$B$ to the class of~$\E(X)$,
hence these two classes have the same order.  Now there
is a natural morphism $r \colon M^1(U) \rightarrow Y$, and
the class of~$\EM(U)$ is the image of the class of~$B$
by the map $\Ext^2(r,\Gm) \colon \Ext^2_{\Sh(\Sm/k)}(Y,\Gm)\rightarrow \Ext^2_{\Sh(\Sm/k)}(M^1(U),\Gm)$ it induces.
The first assertion of Proposition~\ref{propcompareemuex} is thus proved;
moreover, to establish the second assertion, it suffices to show that $\Ext^2(r,\Gm)$ is bijective
as soon as $\NS(U \otimes_k \kbar)=0$.  
Let $\NS_{X/k} \in \Sh(\Sm/k)$ denote the quotient $\Pic_{X/k}/\Pic^0_{X/k}$.
If $\NS(U \otimes_k \kbar)=0$, the natural morphism of étale sheaves $\Div_{X/k,F} \rightarrow \NS_{X/k}$ is surjective.
It follows that~$r$ is a quasi-isomorphism, and hence indeed that $\Ext^2(r,\Gm)$ is bijective.
\end{proof}

\subsection{Main construction}
\label{subsecmaincons}

Before describing the construction alluded to in the introduction of~§\ref{secquat},
we briefly sketch, in the particular case where~$U=X$ and~$X$ is a torsor under an abelian variety~$A$ over~$k$,
a construction which exploits the same idea while avoiding most of the technical details.
Here the goal is to reconstruct~$X$ from~$\EM(X)$ alone; note that the $2$\nobreakdash-extension~$\EM(X)$ is
simply the pullback of~$E(X)$ by the natural map $\Pic^0_{X/k} \rightarrow \Pic_{X/k}$ since~$X$ is proper.  Let $S \in \Sm/k$.  Let $\catExt_S$ denote the
category whose objects are the $2$\nobreakdash-extensions of~$\Pic^0_{X/k}\times_k S$ by~$\GmS$ in $\Sh(\Sm/S)$, and whose morphisms
are the morphisms of complexes which induce the identity on~$\Pic^0_{X/k}\times_k S$ and on~$\GmS$.  Since this is only a sketch, we ignore
set-theoretic problems and pretend that~$\catExt_S$ is a small category (one should really fix a large enough cardinal~$\kappa$ and work
only with sheaves which take values in abelian groups of cardinality~\mbox{$\leq \kappa$}).  We can then consider the geometric realisation~$R_S$ of the nerve of~$\catExt_S$.
This is a topological space.  Objects of~$\catExt_S$ define points of~$R_S$, and connected components of~$R_S$ correspond to Yoneda equivalence classes
of $2$\nobreakdash-extensions.
Let $\pi(S)$ denote the set of homotopy classes of paths on~$R_S$ from the trivial $2$\nobreakdash-extension to~$\EM(X)|_S$.
Note that $\pi(S)\neq\emptyset$ if and only if the class of~$\EM(X)|_S$ in \mbox{$\Ext^2_{\Sh(\Sm/S)}(\Pic^0_{X/k}\times_k S,\GmS)$} is~$0$; in a sense, $\pi(S)$ is the set
of ``Yoneda trivialisations'' of the $2$\nobreakdash-extension $\EM(X)|_S$ up to a certain homotopy relation.
Clearly~$\pi(S)$ depends functorially on~$S$, and thus defines a presheaf in sets on~$\Sm/k$.  Let~$\pitilde$ be the associated étale sheaf.
Then we claim that~$\pitilde$ and~$X$ are canonically isomorphic as étale sheaves on~$\Sm/k$.  To prove it one first reduces to establishing that
if~$X=A$, a certain canonical morphism $\alphatilde \colon A \rightarrow \pitilde$ is an isomorphism. One then
identifies~$\pi(S)$ with the fundamental group of~$R_S$, where the trivial $2$\nobreakdash-extension is chosen as base point.
According to Retakh~\cite[Theorem~1]{retakh}, this fundamental group is canonically isomorphic to \mbox{$\Ext^1_{\Sh(\Sm/S)}(\Pic^0_{X/k}\times_k S, \GmS)$}, which in turn is
canonically isomorphic to~$A(S)$ by the Barsotti-Weil formula (note that $\Pic^0_{X/k}=\Pic^0_{A/k}$).
It only remains to be shown that the resulting canonical isomorphism $A(S) \isoto \pi(S)$ coincides
with~$\alphatilde$, and this is a consequence of the independent results of §\ref{subsecexplicit}.

This construction relies on the description of~$\Ext^n$ groups in abelian categories in terms of $n$\nobreakdash-extensions and Yoneda equivalence.
To carry it out in the general case, one would need a similar description for the hyperext group $\Ext^2_{\Sh(\Sm/S)}(M^1(U)|_S,\GmS)$
as well as a version of Retakh's result which could be applied in this context.  As it turns out, both of these tools are contained in
the paper of Neeman and Retakh~\cite{neemanretakh}; however, making explicit
the construction given in~\cite{neemanretakh} of a canonical isomorphism between the fundamental group
of the nerve of a category analogous to~$\catExt_S$ and the hyperext group $\Ext^1_{\Sh(\Sm/S)}(M^1(U)|_S,\GmS)$, and then
using it to prove that~$\alphatilde$ is an isomorphism, is an excessively daunting task.
Fortunately, a shortcut using derived categories is available, as was explained to the author by Shoham Shamir (see Lemma~\ref{lemmashoham}).
In the proofs below, we thus do not have recourse to the results of~\cite{retakh} or~\cite{neemanretakh}, at the expense of employing a more abstract definition for
the sheaf~$\pitilde$.

\bigskip
\begin{construction}
\label{theconstruction}
Let~$G$ be a semi-abelian variety over~$k$ (so that~$M_1(G)=G$ and~$M^1(G)$ is the $1$\nobreakdash-motive dual to~$G$).
To every $2$\nobreakdash-extension
\begin{equation}
\label{constext}
\xymatrix{
0 \ar[r] & \Gm \ar[r] & E \ar[r] & F \ar[r] & M^1(G) \ar[r] & 0
}
\end{equation}
of~$M^1(G)$ by~$\Gm$ in $C(\Sh(\Sm/k))$,
where~$E$ is concentrated in degree~$0$ and~$F$ in degrees~$-1$ and~$0$, we are going to associate
a $k$\nobreakdash-torsor under~$G$.
\end{construction}

\bigskip
The conditions on~$E$ and~$F$ are not necessary to carry out the construction,
but they make some of the arguments simpler and are unrestrictive enough for our purposes.

Let us first define a category $\catExt_S(G)$ for any $S \in \Sm/k$.  Its objects are the $2$\nobreakdash-extensions
\begin{equation}
\label{eqobjcat}
\xymatrix{
0 \ar[r] & \GmS \ar[r] & E \ar[r]^f & F \ar[r] & M^1(G)|_S \ar[r] & 0
}
\end{equation}
in $C(\Sh(\Sm/S))$ such that~$E$ is concentrated in degree~$0$ and~$F$ in degrees~$-1$ and~$0$.
We define a morphism from~(\ref{eqobjcat}) to
\begin{equation}
\label{eqobjcatp}
\xymatrix{
0 \ar[r] & \GmS \ar[r] & E' \ar[r]^{f'} & F' \ar[r] & M^1(G)|_S \ar[r] & 0
}
\end{equation}
in the category~$\catExt_S(G)$ to be a morphism $C(f) \rightarrow C(f')$ in $D^b(\Sh(\Sm/S))$ such that the diagram
\begin{equation}
\label{eqdiagcomm}
\begin{aligned}
\xymatrix{
\GmS[1] \ar@{=}[d] \ar[r] & C(f) \ar[d] \ar[r] & M^1(G)|_S \ar@{=}[d] \ar[r] & \GmS[2] \ar@{=}[d] \\
\GmS[1] \ar[r] & C(f') \ar[r] & M^1(G)|_S \ar[r] & \GmS[2]
}
\end{aligned}
\end{equation}
commutes, where the rows are the distinguished triangles~(\ref{eqdisttri}) associated with~(\ref{eqobjcat}) and~(\ref{eqobjcatp}).

We note that all morphisms in~$\catExt_S(G)$ are isomorphisms, by the triangulated five lemma (\cite[Cor.~1.5.5]{kascha}).
Furthermore, we note that two objects of~$\catExt_S(G)$ determine the same element of the hyperext group
$\Ext^2_{\Sh(\Sm/S)}(M^1(G)|_S,\GmS)$ if and only if there exists a morphism between them in~$\catExt_S(G)$.

\bigskip
For any $u,v\in\catExt_{\Spec(k)}(G)$ and any $S\in\Sm/k$, we let $\pi_{u,v}(S)=\Hom_{\catExt_S(G)}(u|_S,v|_S)$.
This defines a presheaf in sets~$\pi_{u,v}$ on $\Sm/k$.  Let $\pitilde_{u,v}$ denote the associated étale sheaf.
If $u=v$, then composition defines a group structure on $\pitilde_{u,v}$.

We now consider the case $u=v=\EM(G)$ in more detail.
For any $S \in \Sm/k$ and any $g \in G(S)$, translation by~$g$ on~$G$ defines an automorphism of~$\EM(G)|_S$ in $C(C(\Sh(\Sm/S)))$,
according to Corollary~\ref{coremufunct}.
This automorphism of $\EM(G)|_S$ induces the identity on $M^1(G)|_S$, so that it determines
an element of $\pi_{\EM(G),\EM(G)}(S)$.  Hence we obtain a morphism of presheaves
$\alpha_G \colon G \rightarrow \pi_{\EM(G),\EM(G)}$ on~$\Sm/k$.  We let $\alphatilde_G \colon G \rightarrow \pitilde_{\EM(G),\EM(G)}$ denote the associated morphism of étale sheaves.
Clearly it is a morphism of sheaves in groups.

\bigskip
\begin{theorem}
\label{mainconsthm}
The morphism~$\alphatilde_G$ is an isomorphism.
\end{theorem}

\bigskip
The proof of Theorem~\ref{mainconsthm} is quite involved; we defer it to §\ref{subsecproof}.

We are finally in a position to describe Construction~\ref{theconstruction}.  Let $u \in \A_{\Spec(k)}(G)$.
For any $S \in \Sm/k$, the action of $\pi_{\EM(G),\EM(G)}(S)$
on $\pi_{\EM(G),u}(S)$ by composition is simply transitive.  As a consequence, the sheaf~$\pitilde_{\EM(G),u}$ is a torsor under the sheaf in groups~$\pitilde_{\EM(G),\EM(G)}$.
By Theorem~\ref{mainconsthm}, we can therefore regard it \emph{via}~$\alphatilde_G$ as a torsor under the sheaf defined by~$G$.  Now every étale sheaf on~$\Sm/k$ which is a torsor under~$G$
is representable by a scheme (by Galois descent, see~\cite[6.5/1]{blr});
we thus obtain a $k$\nobreakdash-torsor under~$G$, and this is the torsor that Construction~\ref{theconstruction} associates to~$u$.

We now state and prove the main theorem of this section.

\bigskip
\begin{theorem}
\label{themaintheorem}
Let~$X$ be a smooth proper geometrically integral variety over a field~$k$ and $U \subseteq X$ be a dense open subset.
Then the $k$\nobreakdash-torsor under~$\Alb^0_{U/k}$ associated to~$\EM(U)$ by Construction~\ref{theconstruction}
is canonically isomorphic to~$\Alb^1_{U/k}$.
\end{theorem}

\bigskip
\begin{proof}
The morphism of $1$\nobreakdash-motives $M^1(u_U) \colon M^1(\Alb^1_{U/k}) \rightarrow M^1(U)$ induced by the canonical morphism
$u_U \colon U \rightarrow \Alb^1_{U/k}$ is an isomorphism, since its Cartier dual $M_1(u_U) \colon M_1(U) \rightarrow M_1(\Alb^1_{U/k})$
is clearly an isomorphism.  Therefore it follows from Proposition~\ref{propyondeq} that~$\EM(U)$ and~$\EM(\Alb^1_{U/k})$ are canonically isomorphic
in $\catExt_{\Spec(k)}(\Alb^0_{U/k})$, and hence that the $k$\nobreakdash-torsors associated to
them by Construction~\ref{theconstruction} are canonically isomorphic.
As a consequence, in order to prove Theorem~\ref{themaintheorem}, we may assume $U=\Alb^1_{U/k}$.
Let $G=\Alb^0_{U/k}$.
For any $S \in \Sm/k$,
an element of~$U(S)$ canonically determines an $S$\nobreakdash-isomorphism $U \times_k S \isoto G \times_k S$, hence, by Corollary~\ref{coremufunct},
an isomorphism $\EM(G)|_S \isoto \EM(U)|_S$ in $C(C(\Sh(\Sm/S)))$, hence an element of $\pi_{\EM(G),\EM(U)}(S)$.
We thus obtain a morphism of étale sheaves of sets $\alphatilde_U \colon U \rightarrow \pitilde_{\EM(G),\EM(U)}$.  It only remains to be shown that~$\alphatilde_U$
is an isomorphism.  This question is local for the étale topology on~$\Sm/k$, so that by extending scalars and choosing a point of~$U$, we may
assume that~$U=G$.  In this case $\alphatilde_U=-\alphatilde_G$; applying Theorem~\ref{mainconsthm} concludes the proof.
\end{proof}

\bigskip
\begin{corollary}
\label{quatcor}
Let~$X$ be a smooth proper geometrically integral variety over a field~$k$.
For any dense open subset $U \subseteq X$, the order of the class of~$\EM(U)$ in the hyperext group $\Ext^2_{\Sh(\Sm/k)}(M^1(U),\Gm)$
is equal to~$P(U)$.
The order of the class of~$\E(X)$ in $\Ext^2_{\Sh(\Sm/k)}(\Pic_{X/k},\Gm)$
is equal to~$\Pgen(X)$.
\end{corollary}

\bigskip
\begin{proof}
The second assertion of the corollary follows from the first and from Proposition~\ref{propcompareemuex}, so we need only prove
the first assertion. We fix a dense open subset $U \subseteq X$ and an integer $n \geq 1$,
and let $f \colon U \rightarrow \Alb^n_{U/k}$ denote the composition of
\mbox{$u_U \colon U \rightarrow \Alb^1_{U/k}$} with the canonical morphism $\Alb^1_{U/k} \rightarrow \Alb^n_{U/k}$.
Then $M_1(f) \colon M_1(U) \rightarrow M_1(\Alb^n_{U/k})$
identifies with multiplication by~$n$ on~$M_1(U)$ \emph{via} the canonical isomorphism $M_1(\Alb^n_{U/k})=M_1(U)$;
therefore the map $M^1(f) \colon M^1(\Alb^n_{U/k}) \rightarrow M^1(U)$ identifies with multiplication by~$n$
on~$M^1(U)$. From this we deduce that the class of~$\EM(U)$ in $\Ext^2_{\Sh(\Sm/k)}(M^1(U),\Gm)$
is killed by~$n$ if and only if the class of~$\EM(\Alb^n_{U/k})$ is trivial, thanks to Proposition~\ref{propyondeq}.
The corollary will now be established if we prove that the class of~$\EM(\Alb^n_{U/k})$ is trivial if and only
if $\Alb^n_{U/k}(k)\neq\emptyset$.  One implication is easy: if $\Alb^n_{U/k}(k)\neq\emptyset$, then
Proposition~\ref{propyondeq} applied to a $k$\nobreakdash-morphism $\Spec(k) \rightarrow \Alb^n_{U/k}$ shows
that the class of~$\EM(\Alb^n_{U/k})$ is trivial.  Let us now assume, conversely, that the class of $\EM(\Alb^n_{U/k})$ is trivial.
We have $\Alb^0_{U/k}(k)\neq\emptyset$, hence,
again by Proposition~\ref{propyondeq}, the class of $\EM(\Alb^0_{U/k})$ is trivial.  In particular $\EM(\Alb^0_{U/k})$
and $\EM(\Alb^n_{U/k})$ define the same class in the hyperext group $\Ext^2_{\Sh(\Sm/k)}(M^1(U),\Gm)$.  It follows
that $\pi_{\EM(\Alb^0_{U/k}),\EM(\Alb^n_{U/k})}(k)\neq\emptyset$, which in turn implies
that the $k$\nobreakdash-torsor associated to~$\EM(\Alb^n_{U/k})$ by Construction~\ref{theconstruction}
has a rational point.  By Theorem~\ref{themaintheorem}, this torsor is isomorphic to~$\Alb^n_{U/k}$; hence $\Alb^n_{U/k}(k)\neq\emptyset$,
and the corollary is proved.
\end{proof}

\subsection{Explicit Poincaré sheaves on abelian varieties}
\label{subsecexplicit}

Let~$k$ be a field and~$A$ be an abelian variety over~$k$.  The Poincaré sheaf~$\Prond_A$ is an invertible sheaf on $A \times_k \Pic^0_{A/k}$ which
plays a fundamental rôle in the duality theory of abelian varieties.  One can apply the universal property of the dual abelian variety to~$\Prond_A$ in
two ways, thereby producing a morphism $A \rightarrow \Pic^0_{\Pic^0_{A/k}/k}$ and a morphism $\Pic^0_{A/k} \rightarrow \Pic^0_{A/k}$.
The latter is the identity, by definition of~$\Prond_A$, and the former is the biduality isomorphism \cite[p.~132]{mumfordav}.

Our goal in this subsection is to exhibit an explicit Poincaré sheaf on any abelian variety.  As a consequence we obtain a
very explicit description of the biduality isomorphism.
These results will be used in §\ref{subsecproof} to prove Theorem~\ref{mainconsthm}.

We consider $A \times_k A$ as an $A$\nobreakdash-scheme \emph{via} the second projection.  Let $\Delta \subseteq A \times_k A$ denote the
diagonal and $0 \subseteq A \times_k A$ denote the zero section of the second projection (so that both~$\Delta$ and~$0$ are isomorphic to~$A$
as $k$\nobreakdash-schemes).  Let $V = (A \times_k A) \setminus (\Delta \cup 0)$.
Let
$$d_A \colon \Rstar_{A \times_k A/A,V} \longrightarrow \Div^0_{A \times_k A/A,V} \times \GmA$$
be the morphism of étale sheaves on~$\Sm/A$ defined by $r \mapsto (\mathrm{div}(r), r(0)/r(\Delta))$.

\bigskip
\begin{lemma}
The natural sequence of étale sheaves on~$\Sm/A$
\begin{equation}
\label{poincexplse}
\xymatrix{
0 \ar[r] & \GmA \ar[r] & \Coker(d_A) \ar[r] & \Pic^0_{A \times_k A/A} \ar[r] & 0 \rlap{\text{,}}
}
\end{equation}
where the first map is defined by $x \mapsto (0,x)$ and the second map sends the class of~$(d,x)$ in~$\Coker(d_A)$
to the linear equivalence class of~$d$, is exact.
\end{lemma}

\bigskip
\begin{proof}
Only exactness on the right requires a proof. It suffices
to check that for any $T \in \Sm/A$, any $t \in T$, and any $D \in \Div(A \times_k T)$, there is an open subset~$U$
of~$T$ containing~$t$ such that the restriction of~$D$ to $A \times_k U$ is linearly equivalent, on $A \times_k U$, to a divisor
whose support is contained in $V \times_A U$; and this follows from the vanishing of the Picard group of
the semilocal ring of~$A \times_k T$ at the two points $0 \times t$ and $s(t) \times t$, where $s \colon T \rightarrow A$
denotes the structural morphism of~$T$
(the Picard group of any semilocal ring vanishes).
\end{proof}

\bigskip
The exact sequence~(\ref{poincexplse}) defines a torsor under~$\Gm$ over $\Pic^0_{A \times_k A/A}$
(whose sections over an open subset~$U$ of $\Pic^0_{A \times_k A/A}$ are the elements of $\Coker(d_A)(U)$ which lift
the inclusion $U \subseteq \Pic^0_{A \times_k A/A}$ seen as
an element of $\Pic^0_{A\times_k A /A}(U)$; note that $U \in \Sm/A$),
hence an invertible sheaf~$\Qrond_A$ on $\Pic^0_{A \times_kA/A}=\Pic^0_{A/k}\times_kA$.

To identify line bundles, invertible sheaves, and torsors under~$\Gm$, we use covariant conventions:
the line bundle associated with an invertible sheaf~$\Lrond$ on a scheme~$S$ is the (relative) spectrum
of the symmetric algebra of the dual of~$\Lrond$; the torsor under~$\Gm$ associated with a line bundle $L \rightarrow S$
is the complement of the zero section in~$L$, with the induced action of~$\Gm$.

\bigskip
\begin{theorem}
\label{thpoincare}
The invertible sheaf~$\Qrond_A$ is isomorphic to the Poincaré sheaf~$\Prond_A$.
\end{theorem}

\bigskip
\begin{corollary}
\label{corbiduality}
The biduality isomorphism $A \rightarrow \Pic^0_{\Pic^0_{A/k}/k}$ sends a point $a \in A(k)$ to the class of the invertible
sheaf on~$\Pic^0_{A/k}$ defined by the exact sequence of étale sheaves on~$\Sm/k$
\begin{equation}
\label{poincexplse2}
\xymatrix{
0 \ar[r] & \Gm \ar[r] & \Coker(\Rstar_{A/k,A \setminus\{a,0\}} \rightarrow \Div^0_{A/k,A \setminus\{a,0\}} \times \Gm) \ar[r] & \Pic^0_{A/k} \ar[r] & 0\rlap{\text{,}}
}
\end{equation}
where the map $\Rstar_{A/k,A \setminus\{a,0\}} \rightarrow \Div^0_{A/k,A \setminus\{a,0\}} \times \Gm$ is defined by $r \mapsto (\mathrm{div}(r),r(0)/r(a))$.
\end{corollary}

\bigskip
\begin{proof}[ of Theorem~\ref{thpoincare}]%
For $a=0$, the exact sequence~(\ref{poincexplse2}) is clearly split.  In other words, the morphism $A \rightarrow \Pic_{\Pic^0_{A/k}/k}$
which corresponds to~$\Qrond_A$ by the universal property of the Picard functor maps~$0$ to~$0$.  Hence it induces
a morphism of abelian varieties $A \rightarrow \Pic^0_{\Pic^0_{A/k}/k}$.
Let $\sigma_A \in \End_k(A)$ denote the composition of this morphism
with the inverse of the biduality isomorphism $A \isoto \Pic^0_{\Pic^0_{A/k}/k}$.
The assertion that~$\Prond_A$ and~$\Qrond_A$ are isomorphic now amounts to the equality $\sigma_A=\Id$.

There exist a Jacobian variety~$J$ over~$k$ and a smooth surjective morphism $m \colon J \rightarrow A$ (\cite[Corollary~2.5]{gabberfilling}).
It is straightforward to check that the pullbacks of~$\Qrond_A$ and~$\Qrond_J$ to~$\Pic^0_{A/k} \times_k J$ are isomorphic
(note that $J \in \Sm/A$).  As a consequence, the square
$$
\xymatrix{
J \ar[r]^m \ar[d]^{\sigma_J} & A \ar[d]^{\sigma_A} \\
J \ar[r]^m & A
}
$$
commutes; hence $\sigma_J=\Id$ implies $\sigma_A=\Id$.

We can therefore assume that~$A$ is a Jacobian variety, \emph{i.e.}, that there exists a smooth proper geometrically connected curve~$C$ over~$k$
such that $A=\Pic^0_{C/k}$.

After replacing~$k$ with a finite separable extension (which, for the purpose of proving the equality of endomorphisms $\sigma_A=\Id$, is harmless),
we may assume that $C(k)\neq\emptyset$.  Fix a point $c_0 \in C(k)$.
There is a canonical embedding $\iota \colon C \rightarrow A$ such that $\iota(c_0)=0$ (\cite[§2]{milnejac}).
In order to prove that $\sigma_A = \Id$, it is enough to prove that $\sigma_A \circ \iota = \iota$, since the image of~$\iota$ generates~$A$.
As~$C$ is reduced, it even suffices to establish that $\sigma_A(\iota(c))=\iota(c)$
for every closed point $c \in C$ whose residue field is separable over~$k$,
or equivalently that $\Prond_A|_{\Pic^0_{A/k} \times \{\iota(c)\}}$
and
$\Qrond_A|_{\Pic^0_{A/k} \times \{\iota(c)\}}$
are isomorphic invertible sheaves for every such $c \in C$.  Fix a closed point $c \in C$ with separable residue field over~$k$.
After replacing~$k$ with the residue field of~$c$, we may assume that $c \in C(k)$.
The morphism $\Pic^0_{A/k} \rightarrow \Pic^0_{C/k}=A$ induced by~$\iota$ is an isomorphism (\cite[§6]{milnejac}).
Let $\iota' \colon C \rightarrow \Pic^0_{A/k}$ denote the composition of~$\iota$ with the inverse of this isomorphism.
Again by~\cite[§6]{milnejac}, two algebraically equivalent invertible sheaves on $\Pic^0_{A/k}$ are isomorphic if and only if their inverse
images by~$\iota'$ are isomorphic.  In particular, we need only prove
that $s^\star \Prond_A$ and $s^\star \Qrond_A$ are isomorphic,
where $s \colon C \rightarrow \Pic^0_{A/k} \times_k A$ is the morphism defined by $s(x)=(\iota'(x),\iota(c))$.

According to \cite[§6, 6.11]{milnejac}, the invertible sheaf $s^\star \Prond_A$ is isomorphic to $\Orond(c-c_0)$.
Let $X \rightarrow C$ and $Y \rightarrow C$ denote the torsors under~$\Gm$ respectively
defined by the invertible sheaves~$\Orond(c-c_0)$ and $s^\star \Qrond_A$.
The embedding~$\iota$ induces a morphism from the exact sequence~(\ref{poincexplse2}) with $a=\iota(c)$
to the exact sequence
$$
\xymatrix{
0 \ar[r] & \Gm \ar[r] & \Coker\left(\Rstar_{C/k,C \setminus\{c,c_0\}} \xrightarrow{\;\;d_C\;} \Div^0_{C/k,C \setminus\{c,c_0\}} \times \Gm\right) \ar[r] & \Pic^0_{C/k} \ar[r] & 0\rlap{\text{,}}
}
$$
where~$d_C$ is defined by $r \mapsto (\mathrm{div}(r),r(c_0)/r(c))$.
This is an isomorphism of exact sequences, by the five lemma.  Hence the torsor $Y \rightarrow C$ can be described as follows: its sections
over an open subset $U \subseteq C$ are the liftings in $\Coker(d_C)(U)$ of the element
of~$\Pic^0_{C/k}(U)$ defined by the composition
\mbox{$U \subseteq C \xrightarrow{\;\!\iota} A=\Pic^0_{C/k}$}
(note that $U \in \Sm/k$).
The torsor $X \rightarrow C$ is easier to describe: its sections over an open subset $U \subseteq C$ are the rational functions $g \in k(C)^\star$
such that the divisor $\mathrm{div}(g)+c-c_0$ is supported on~$C \setminus U$.

We shall now exhibit a $\Gm$\nobreakdash-equivariant $C$\nobreakdash-morphism $\phi \colon X \rightarrow Y$.
Let us fix once and for all a rational function $\pi \in k(C)^\star$ which is a uniformiser at~$c_0$ and which moreover is invertible at~$c$
if $c\neq c_0$. (The existence of~$\pi$ follows from the vanishing of the Picard group of the semilocal ring of~$C$ at~$c$ and~$c_0$.)
We define~$\phi$ in the neighbourhood of an arbitrary closed point $c_1 \in C$.  Let $U \subseteq C$ be an open subset containing~$c_1$.
Let $D \in \Div(C \times_k U)$ be a divisor linearly equivalent to $\Delta - (c_0 \times U)$ (where~$\Delta$ denotes the image of the diagonal
embedding $U \hookrightarrow C \times_k U$), with disjoint support from $(c \times c_1) \cup (c_0 \times c_1)$.
Such a divisor exists because the Picard group of a semilocal ring is trivial.  By shrinking~$U$ if necessary, we may assume that~$D$ has support
disjoint from $(c \times U) \cup (c_0 \times U)$.  Let $h \in k(C \times_k U)^\star$ satisfy $\mathrm{div}(h)=\Delta - (c_0 \times U) - D$, and
let $h' \in k(C \times_k U)^\star$ be the rational function defined by $h'(x,y)=h(x,y)\pi(x)$.
The divisor of~$h'$ coincides with~$\Delta$ in a neighbourhood of $(c \times U) \cup (c_0 \times U)$.  As a consequence,
the rational function $y \mapsto g(y)h'(c,y)/h'(c_0,y)$ is invertible on~$U$
for any $g \in k(C)^\star$ such that $\mathrm{div}(g)+c-c_0$ is supported on~$C \setminus U$.
We let $\phi(U) \colon X(U) \rightarrow Y(U)$ be the map which sends~$g$ to the class in $\Coker(d_C)(U)$ of the pair consisting
of the divisor $D \in \Div(C \times_k U)$ and the function $(y \mapsto g(y)h'(c,y)/h'(c_0,y)) \in \Gm(U)$.
Note that the restrictions of~$D$ to the geometric fibres of the second projection $C \times_k U \rightarrow U$ are indeed algebraically equivalent to~$0$
and supported on $C \setminus \{c,c_0\}$.

Suppose~$D_1 \in \Div(C \times_k U)$ and~$h_1 \in k(C \times_k U)^\star$ are another divisor and another rational function satisfying the same
requirements as~$D$ and~$h$, \emph{i.e.},
$D_1$ is linearly equivalent to
$\Delta - (c_0 \times U)$, has support disjoint from $(c \times U) \cup (c_0 \times U)$, and satisfies
$\mathrm{div}(h_1)=\Delta - (c_0 \times U) - D_1$.   Then,
for any $g \in k(C)^\star$ such that $\mathrm{div}(g)+c-c_0$ is supported on~$C \setminus U$, the pairs $(D, gh'(c,-)/h'(c_0,-))$ and $(D_1, gh'_1(c,-)/h'_1(c_0,-))$,
where $h'_1(x,y)=h_1(x,y)\pi(x)$, define the same class in $\Coker(d_C)(U)$.  Indeed, letting $r=h/h_1$, these two pairs
differ by $(\mathrm{div}(r), r(c_0,-)/r(c,-))$, which
by definition of~$d_C$ vanishes in $\Coker(d_C)(U)$.
Hence the map~$\phi(U)$ does not depend on the choices of~$D$ and~$h$.  In addition, it is clearly $\Gm(U)$\nobreakdash-equivariant
and functorial with respect to~$U$.  We have therefore defined a $\Gm$\nobreakdash-equivariant morphism $\phi \colon X \rightarrow Y$
in the category of Zariski sheaves of sets on~$C$.
Now by Grothendieck's Hilbert~90 theorem, giving a $\Gm$\nobreakdash-equivariant $C$\nobreakdash-morphism between
torsors under~$\Gm$ over~$C$ is equivalent to giving a $\Gm$\nobreakdash-equivariant morphism between the Zariski sheaves on~$C$ they define
(see~\cite[Arcata, II-4]{ega4h}).
Moreover, such a morphism is necessarily an isomorphism; hence the proof of Theorem~\ref{thpoincare} is complete.
\end{proof}

\subsection{Proof of Theorem~\ref{mainconsthm}}
\label{subsecproof}

We now show that the morphism of sheaves $\alphatilde_G$ is an isomorphism.
This question is local for the étale topology on~$\Sm/k$.  By replacing~$k$ with a finite separable extension, we can therefore
assume (for simplicity) that~$G$ is an extension of an abelian variety by a split torus.

The following lemma and its proof are a variant for the hyperext groups under consideration of a general result on
$\Ext$ groups in abelian categories which was communicated to the author by Shoham Shamir.

\bigskip
\begin{lemma}
\label{lemmashoham}
Fix $S \in \Sm/k$ and $u \in \catExt_S(G)$.
For $\phi \in \Hom_{D^b(\Sh(\Sm/S))}(M^1(G)|_S,\GmS[1])$, we let~$A(\phi)$ denote the composition
$C(f) \rightarrow M^1(G)|_S \xrightarrow{\phi} \GmS[1] \rightarrow C(f)$, where the unlabeled maps
are taken from the exact triangle~(\ref{eqdisttri}) associated with~$u$.  Then $\phi \mapsto \Id + A(\phi)$ defines an isomorphism
from the hyperext group $\Ext^1_{\Sh(\Sm/S)}(M^1(G)|_S,\GmS)$ to the group of automorphisms of~$u$
in~$\catExt_S(G)$.
\end{lemma}

\bigskip
\begin{proof}
The subscript in $\Hom_{D^b(\Sh(\Sm/S))}(-,-)$ will be dropped for the rest of the proof.
Let us denote by $\rho \colon \GmS[1] \rightarrow C(f)$ and $\sigma \colon C(f) \rightarrow M^1(G)|_S$ the maps
which appear in the exact triangle~(\ref{eqdisttri}) associated with~$u$, so that $A(\phi)=\rho \circ \phi \circ \sigma$.
First it is immediate that $\Id+A(\phi)$ is an endomorphism of~$u$ in~$\catExt_S(G)$, and hence an automorphism.
We have $A(\phi) \circ A(\psi)=0$ and $A(\phi+\psi)=A(\phi)+A(\psi)$ for any~$\phi$ and~$\psi$, therefore
$\phi \mapsto \Id + A(\phi)$ is a group morphism.
We note that $\Hom(\GmS[2], C(f))=0$ and $\Hom(M^1(G)|_S, M^1(G)|_S[-1])=0$.
Indeed the first group vanishes because~$C(f)$ is concentrated in degrees~$-1$ and~$0$ (quite generally one has $\Hom(E,F)=0$
if~$E$ is concentrated in degrees~$\leq n$ and~$F$ in degrees $\geq n+1$, for some~$n \in \Z$);
similarly, $M^1(G)|_S$ being concentrated in degrees~$-1$ and~$0$ implies that the second group injects into the group of homomorphisms
from an $S$\nobreakdash-abelian scheme to an $S$\nobreakdash-lattice, hence is trivial.
Taking these remarks into account while applying the functors $\Hom(-,C(f))$ and $\Hom(M^1(G)|_S,-)$ to the distinguished triangle~(\ref{eqdisttri}),
we obtain the two exact sequences
\begin{equation}
\label{shohames1}
\xymatrix{
0 \ar[r] & \Hom(M^1(G)|_S,C(f)) \ar[r]^(.53){-\circ\sigma} & \Hom(C(f),C(f)) \ar[r]^(.465){-\circ\rho} & \Hom(\GmS[1], C(f))
}
\end{equation}
and
\begin{equation}
\label{shohames2}
\xymatrix{
0 \ar[r] & \Hom(M^1(G)|_S,\GmS[1]) \ar[r]^(.53){\rho\circ-} & \Hom(M^1(G)|_S, C(f)) \ar[r]^(.465){\sigma\circ-} & \Hom(M^1(G)|_S, M^1(G)|_S)\rlap{\text{.}}
}
\end{equation}
It follows that $\phi \mapsto A(\phi)$ induces a bijection from $\Hom(M^1(G)|_S,\GmS[1])$
to the group of morphisms $a \colon C(f) \rightarrow C(f)$ such that $a \circ \rho =0$ and $\sigma\circ a=0$.
Now a morphism $a \colon C(f) \rightarrow C(f)$ satisfies these conditions if and only if $\Id + a$
is an automorphism of~$u$ in $\catExt_S(G)$, hence the lemma.
\end{proof}

\bigskip
Let~$X$ denote a smooth compactification of~$G$ and $F=X \setminus G$.
For $S \in \Sm/k$, we denote by $f \colon \Rstar_{X \times_k S/S} \rightarrow [\Div^0_{X \times_k S/S,F \times_k S} \rightarrow \Div^0_{X \times_k S/S}]$ the restriction to~$\Sm/S$ of the morphism
of complexes of étale sheaves on~$\Sm/k$ which appears in the definition of~$\EM(G)$ (see §\ref{subsecmotivic}).
Let $$\beta_G \colon \Ext^1_{\Sh(\Sm/S)}(M^1(G)|_S,\GmS) \longisoto
\pi_{\EM(G),\EM(G)}(S)$$ be the isomorphism given by Lemma~\ref{lemmashoham}.  By abuse of notation we shall write~$\alpha_G$
for the map $\alpha_G(S) \colon G(S) \rightarrow \pi_{\EM(G),\EM(G)}(S)$ induced by~$\alpha_G$.

\bigskip
\begin{lemma}
\label{lemmaabfunct}
Let $S \in \Sm/k$.
Let~$G'$ be another semi-abelian variety over~$k$ and $m \colon G \rightarrow G'$ be a homomorphism.
The
square
$$
\xymatrix@R=8ex@C=23ex{
G(S) \ar[r]^m \ar[d]^{\beta_G^{-1} \circ \alpha_G} & G'(S) \ar[d]^{\beta_{G'}^{-1} \circ \alpha_{G'}} \\
\Ext^1_{\Sh(\Sm/S)}(M^1(G)|_S,\GmS) \ar[r]^{\Ext^1_{\Sh(\Sm/S)}(M^1(m),\GmS)} & \Ext^1_{\Sh(\Sm/S)}(M^1(G')|_S,\GmS)
}
$$
commutes.
\end{lemma}

\bigskip
\begin{proof}
We fix a smooth compactification~$X'$ of~$G'$ and
let $$f' \colon \Rstar_{X' \times_k S/S} \longrightarrow [\Div^0_{X' \times_k S/S,F' \times_k S} \rightarrow \Div^0_{X' \times_k S/S}]$$
be the restriction to~$\Sm/S$ of the morphism of complexes which appears in the definition of~$\EM(G')$, where $F' = X'\setminus G'$.
Let $g \in G(S)$.
Translation by~$g$ (resp.~by~$m(g)$) on~$G$ (resp.~on~$G'$) induces an automorphism of the mapping cone~$C(f)$ (resp.~$C(f')$)
in $D^b(\Sh(\Sm/S))$, which we denote~$t^\star_g$ (resp.~$t^\star_{m(g)}$).
Let $\rho \colon \GmS[1] \rightarrow C(f)$, $\rho' \colon \GmS[1] \rightarrow C(f')$,
$\sigma \colon C(f) \rightarrow M^1(G)|_S$, and~$\sigma' \colon C(f') \rightarrow M^1(G')|_S$
denote the maps defined by the exact triangles~(\ref{eqdisttri}) associated with~$\EM(G)|_S$ and~$\EM(G')|_S$.

By Lemma~\ref{lemmashoham}, there exist unique morphisms $\phi \colon M^1(G)|_S \rightarrow \GmS[1]$
and~$\phi' \colon M^1(G')|_S \rightarrow \GmS[1]$ in $D^b(\Sh(\Sm/S))$ such that $\rho \circ \phi \circ \sigma = t^\star_g-\Id$
and $\rho'\circ\phi'\circ\sigma'=t^\star_{m(g)}-\Id$.
We denote by $m^\star \colon C(f') \rightarrow C(f)$ and $M^1(m) \colon M^1(G')|_S \rightarrow M^1(G)|_S$ the morphisms in $D^b(\Sh(\Sm/S))$ induced by~$m$
(the first one being given by Proposition~\ref{propyondeq}).
The conclusion of Lemma~\ref{lemmaabfunct} amounts to the equality \mbox{$\phi \circ M^1(m)=\phi'$} in $\Hom_{D^b(\Sh(\Sm/S))}(M^1(G')|_S,\GmS[1])$.
Clearly one has $t_g^\star \circ m^\star = m^\star \circ t^\star_{m(g)}$,
$m^\star \circ \rho' = \rho$, and $M^1(m) \circ \sigma' = \sigma \circ m^\star$.  Hence $\rho \circ \phi \circ M^1(m) \circ \sigma' = \rho \circ \phi' \circ \sigma'$.
We now argue as in the proof of Lemma~\ref{lemmashoham} to deduce from this equality
first that $\rho \circ \phi \circ M^1(m)=\rho \circ \phi'$
(using the exact sequence obtained from~(\ref{shohames1}) by replacing the second occurrence of~$C(f)$ with~$C(f')$,
and~$G$, $\sigma$, $\rho$ with~$G'$, $\sigma'$,~$\rho'$),
and then that $\phi \circ M^1(m)=\phi'$
(using the exact sequence obtained from~(\ref{shohames2}) by replacing the first three occurrences of~$G$ with~$G'$).
\end{proof}

\bigskip
The morphism $\beta_G^{-1} \circ \alpha_G$ depends functorially on~$S$, and can therefore be regarded as a morphism of presheaves.
By considering the associated morphism of sheaves,
one deduces from Lemma~\ref{lemmaabfunct} and from the five lemma that in order to prove that $\alphatilde_G$ is an isomorphism,
it is harmless to assume that~$G$ is either a split torus or an abelian variety, and then that~$G$ is either~$\Gm$ or an abelian variety (and accordingly $X=\P^1_k$ or $X=G$).
Under this assumption, we shall establish that the morphism of presheaves~$\alpha_G$ is itself an isomorphism.

We first suppose~$G$ is an abelian variety, so that $M^1(G)=\Pic^0_{G/k}$. For any $S \in \Sm/k$, the Barsotti-Weil morphism
$$\Ext^1_{\Sh(\Sm/S)}(\Pic^0_{G/k} \times_k S, \GmS) \longrightarrow \Pic^0_{\Pic^0_{G/k}/k}(S)\rlap{\text{,}}$$
which sends an extension of~$\Pic^0_{G/k} \times_k S$ by~$\GmS$ to the class of the torsor under~$\GmS$ over~$\Pic^0_{G/k} \times_k S$ that it defines,
is an isomorphism (see~\cite[III.18]{oort}).  Composing its inverse with
the biduality isomorphism $G(S) \isoto \Pic^0_{\Pic^0_{G/k}/k}(S)$ yields
an isomorphism $\gamma_G \colon G(S) \isoto \Ext^1_{\Sh(\Sm/S)}(M^1(G)|_S,\GmS)$.
Clearly $\gamma_G^{-1} \circ \beta_G^{-1} \circ \alpha_G$ depends functorially on~$S$, hence defines a morphism
of presheaves $G \rightarrow G$ on~$\Sm/k$.  By Yoneda's lemma, this morphism comes from an endomorphism~$m$ of
the abelian variety~$G$.  For the morphism of presheaves~$\alpha_G$ to be an isomorphism, it suffices that $m=\Id$.
Now two endomorphisms of a smooth algebraic $k$\nobreakdash-group coincide as soon as they coincide
on $\ell$\nobreakdash-points for all finite separable extensions~$\ell/k$.  As a consequence, we will be done if
the morphisms of groups $G(\ell) \rightarrow \pi_{\EM(G),\EM(G)}(\ell)$ induced
by~$\alpha_G$ and by~$\beta_G \circ \gamma_G$ are equal for all finite separable extensions $\ell/k$; this is what we shall
prove now.  After extending the scalars from~$k$ to~$\ell$, we may assume $\ell=k$.

Let $g \in G(k)$.  Let
$d \colon \Rstar_{G/k,G\setminus \{g,0\}} \rightarrow \Div^0_{G/k,G\setminus \{g,0\}} \times \Gm$
denote the morphism of étale sheaves on~$\Sm/k$ defined by $r \mapsto (\mathrm{div}(r),r(0)/r(g))$.
By Corollary~\ref{corbiduality}, the class in
$\Ext^1_{\Sh(\Sm/k)}(M^1(G),\Gm)$ of the natural exact sequence
\begin{equation}
\label{exsegammag}
\xymatrix{
0 \ar[r] & \Gm \ar[r] & \Coker(d) \ar[r] & \Pic^0_{G/k} \ar[r] & 0
}
\end{equation}
is equal to $\gamma_G(g)$.  Let $\phi \colon \Pic^0_{G/k} \rightarrow \Gm[1]$ be the morphism in $D^b(\Sh(\Sm/k))$ corresponding
to~(\ref{exsegammag}).  We keep the notations~$A$, $\rho$, and~$\sigma$ from the statement and the proof of Lemma~\ref{lemmashoham} with $S=\Spec(k)$.
It only remains to be shown that $\Id + A(\phi)$ coincides with the morphism $t^\star_g\colon C(f) \rightarrow C(f)$ induced by translation by~$g$
on~$G$.  The complex $C(f)$ is equal to $[\Rstar_{G/k} \rightarrow \Div^0_{G/k}]$ (with~$\Div^0_{G/k}$ in degree~$0$).
The morphism $A(\phi) \in \Hom_{D^b(\Sh(\Sm/k))}(C(f),C(f))$ can be written as a composition of morphisms of complexes and formal inverses
of quasi-isomorphisms of complexes in the following way:
\begin{equation}
\label{aphicomposition}
\begin{aligned}
\xymatrix{
\Rstar_{G/k} \ar[r] \ar[d] & 0 \ar[d] & \Gm \ar[l] \ar[d] \ar[r]^{\Id} & \Gm \ar[r] \ar[d] & \Rstar_{G/k} \ar[d] \\
\Div^0_{G/k} \ar[r] & \Pic^0_{G/k} & \Coker(d) \ar[l] \ar[r] & 0 \ar[r] & \Div^0_{G/k}\rlap{\kern.1em\text{.}}
}
\end{aligned}
\end{equation}
The natural commutative cube
$$
\xymatrix@R=2ex@C=3ex{
& 0 \ar[dd] && \Gm \ar[ll] \ar[dd] \\
\Rstar_{G/k} \ar[ur] \ar[dd] && \Rstar_{G/k,G\setminus \{g,0\}} \ar[dd] \ar[ll] \ar[ur] \\
& \Pic^0_{G/k} && \Coker(d) \ar[ll] \rlap{\text{,}} \\
\Div^0_{G/k} \ar[ur] && \Div^0_{G/k,G\setminus \{g,0\}} \ar[ll] \ar[ur]
}
$$
where the map $\Rstar_{G/k,G\setminus \{g,0\}} \rightarrow \Gm$ is defined by $r \mapsto r(g)/r(0)$, enables one to
rewrite~(\ref{aphicomposition}) as the composition of a single formal inverse of a quasi-isomorphism followed by a single morphism of complexes.
Indeed the frontmost and the backmost faces, regarded as leftwards morphisms of vertical complexes, are quasi-isomorphisms.
From this explicit description of~$A(\phi)$, one deduces that the morphism
$\Id + A(\phi) - t^\star_g$ is equal to~$0$
in $D^b(\Sh(\Sm/k))$ if and only if the morphism of (vertical) complexes
\begin{equation}
\label{mornullhom}
\begin{aligned}
\xymatrix@C=8ex{
\Rstar_{G/k,G\setminus \{g,0\}} \ar[d] \ar[r] & \Rstar_{G/k} \ar[d] \\
\Div^0_{G/k,G\setminus \{g,0\}} \ar[r]^(0.6){\Id - t^\star_g} & \Div^0_{G/k}\rlap{\text{,}}
}
\end{aligned}
\end{equation}
where the top horizontal map sends a rational function~$r$ to the rational function $$x \longmapsto \frac{r(g)r(x)}{r(0)r(x+g)}\rlap{\text{,}}$$
is equal to~$0$ in $D^b(\Sh(\Sm/k))$.
For any $T \in \Sm/k$ and any $D \in \Div(G \times_k T)$
whose support is contained in $(G\setminus \{g,0\}) \times_k T$ and whose restriction to every geometric fibre of the second projection $G \times_k T \rightarrow T$
is algebraically equivalent to~$0$, the section of $\Pic^0_{G/k}$ on~$T$ defined by~$D-t^\star_g D$ is trivial.  Indeed, as is well known,
the linear equivalence class of a divisor on an abelian variety is invariant under translations if the divisor is algebraically equivalent to~$0$
(see \cite[p.~75]{mumfordav}).  Therefore there exists (locally on~$T$ for the Zariski topology) a rational function~$r$ on~$G \times_k T$
whose divisor is $D-t^\star_g D$.  Let~$h(D)$ be the unique such~$r$ which in addition satisfies~$r(0)=1$.
The association of~$h(D)$ to~$D$ defines a morphism of étale sheaves $h \colon \Div^0_{G/k,G\setminus \{g,0\}} \rightarrow \Rstar_{G/k}$ on~$\Sm/k$.
The diagram obtained by adjoining~$h$ to the square~(\ref{mornullhom}) is still commutative.
As a consequence, the morphism of complexes~(\ref{mornullhom}) is null-homotopic, hence is indeed equal to~$0$ in $D^b(\Sh(\Sm/k))$.

The case where $G=\Gm$ and $X=\P^1_k$ is similar, though substantially easier; we explain it briefly.
Let $S \in \Sm/k$.
We have $M^1(G)|_S=\Z[1]$ in $C(\Sh(\Sm/S))$, where~$\Z$ denotes the constant sheaf~$\Z$ on $\Sm/S$.
In order to prove that~$\alpha_G$ is an isomorphism, it suffices to check that for any $g \in G(S)$,
the equality $t^\star_g = \Id+A(\phi)$ holds in $\Hom_{D^b(\Sh(\Sm/S))}(C(f),C(f))$, where $\phi \colon M^1(G)|_S \rightarrow \GmS[1]$
now denotes the morphism in $D^b(\Sh(\Sm/S))$ obtained by shifting the morphism of sheaves $\Z \rightarrow \GmS$ which maps~$1$ to~$g$.
For this it suffices that the morphism of (vertical) complexes
$$
\xymatrix{
\Z \oplus \Rstar_{\P^1_S/S} \ar[d] \ar[r] & \Z \oplus \Rstar_{\P^1_S/S} \ar[d] \\
\Div^0_{\P^1_S/S} \ar[r]^{\Id - t^\star_g} & \Div^0_{\P^1_S/S}\rlap{\text{,}}
}
$$
where the top horizontal map sends $n \oplus r$ to $0 \oplus (g^n r/t^\star_g r)$
and both vertical maps send $n \oplus r$ to $\mathrm{div}(r) + n(\infty - 0)$,
be null-homotopic.  It is easy to see that a homotopy is given by the morphism $\Div^0_{\P^1_S/S} \rightarrow \Z \oplus \Rstar_{\P^1_S/S}$
defined by $D \mapsto 0 \oplus (r/t^\star_g r)$, where~$r$ is any rational function such that $\mathrm{div}(r)=D$.
Thus the proof of Theorem~\ref{mainconsthm} is complete.

\appendix\section{}
Except for the applications given in~§\ref{sectarith}, most of the present paper deals with smooth varieties over an arbitrary field~$k$.
Care has been taken not to exclude imperfect fields from our treatment, with global fields of positive characteristic in mind as possible choices for~$k$.
To the best of our knowledge, however, the results we use concerning universal morphisms to semi-abelian varieties cannot be found in the literature without the assumption that~$k$ be perfect.
We provide in this appendix arguments to circumvent the few issues which show up when dealing with imperfect fields.

\bigskip
\begin{theorem}[ (Serre)]%
\label{thalbexiste}
Let~$k$ be a field and~$X$ be a geometrically integral variety over~$k$, endowed with a rational point $x \in X(k)$.
There exist a semi-abelian variety~$A$ over~$k$
and a $k$\nobreakdash-morphism $u \colon X \rightarrow A$ mapping~$x$ to~$0$, such that for any semi-abelian
variety~$B$ over~$k$, any $k$\nobreakdash-morphism
$X \rightarrow B$ mapping~$x$ to~$0$ factors uniquely through~$u$.
\end{theorem}

\bigskip
Let~$\kalg$ be an algebraic closure of~$k$.
Theorem~\ref{thalbexiste} is due to Serre when $k=\kalg$ (\cite[Théorème~7]{serremunivalb}).
To prove it in general, we shall follow closely some of the arguments of~\cite{serremunivalb}, and then use the fact that Theorem~\ref{thalbexiste} is true
for~$\kalg$ to shortcut the rest of the proof.  While no condition on~$k$ is imposed in the statement of Theorem~\ref{thalbexiste}, the case of interest is
really that of a separably closed field.

\bigskip
\begin{sketchofproof}
If~$A$ is a semi-abelian variety over~$k$, a $k$\nobreakdash-morphism $a \colon X \rightarrow A$ mapping~$x$ to~$0$
is said to be \emph{generating} if~$A$ is the smallest subgroup scheme of~$A$ over~$k$ through which~$a$ factors.
It is said to be \emph{maximal} if for any semi-abelian variety~$B$ over~$k$ and any $k$\nobreakdash-morphism $b \colon X \rightarrow B$
mapping~$x$ to~$0$, any finite $k$\nobreakdash-homomorphism $h \colon B \rightarrow A$ such that $h \circ b = a$ is an isomorphism.

\bigskip
\begin{lemma}
\label{lemmagenerating}
Let~$A$ be a semi-abelian variety over~$k$ and $a \colon X \rightarrow A$ be a $k$\nobreakdash-morphism mapping~$x$ to~$0$.
If~$a$ is generating, then
$a \otimes_k \kalg \colon X \otimes_k \kalg \rightarrow A \otimes_k \kalg$ is generating.
\end{lemma}

\bigskip
\begin{proof}
Let $a^n \colon X^n \rightarrow A$ be the map $(\uplet{x_1}{x_n})\mapsto \sum a(x_i)$, where~$X^n$ denotes the $n$\nobreakdash-fold
fibred product of~$X$ with itself above~$k$.  The reasoning of~\cite[p.~10-01]{serremunivalb} (in terms of scheme-theoretic images)
shows that~$a$ is generating if and only if~$a^n$ is dominant for large~$n$, and this last condition is clearly preserved
by extension of scalars from~$k$ to~$\kalg$.
\end{proof}

\bigskip
\begin{lemma}
\label{ersatzthun}
Let~$A$ be a semi-abelian variety over~$k$ and $a \colon X \rightarrow A$ be a $k$\nobreakdash-morphism mapping~$x$ to~$0$.
There exist a semi-abelian variety~$B$ over~$k$, a finite $k$\nobreakdash-homomorphism $h \colon B \rightarrow A$, and
a $k$\nobreakdash-morphism $b \colon X \rightarrow B$ mapping~$x$ to~$0$, such that $h \circ b = a$ and~$b$ is maximal.
\end{lemma}

\bigskip
\begin{proof}
Let~$A'$ be the smallest subgroup scheme of~$A$ over~$k$ through which~$a$ factors.  Let~$A''$ denote the reduced closed subscheme of~$A'$
whose underlying topological space is the connected component of~$0$ in~$A'$.  The $k$\nobreakdash-morphism~$a$ factors through~$A''$ since~$X$ is reduced and connected.
Now~$A''$ is a subgroup scheme of~$A'$ over~$k$, according to \cite[Proposition~3.1]{tdte6}.  Hence $A''=A'$, so that~$A'$ is connected
and reduced.  It now follows from \cite[Proposition~3.1]{tdte6} that~$A'$ is smooth over~$k$.  Any smooth connected subgroup scheme of a semi-abelian
variety is itself a semi-abelian variety.  We may therefore replace~$A$ with~$A'$, and assume that~$a$ is generating.  From this point on,
the proof of \cite[Théorème~1]{serremunivalb} applies word for word.
\end{proof}

\bigskip
The arguments in the second part of the proof of \cite[Théorème~2]{serremunivalb} require no modification, except that
the reference to \cite[Théorème~1]{serremunivalb} must be replaced with a reference to Lemma~\ref{ersatzthun} above.
They imply that in order to establish Theorem~\ref{thalbexiste}, it suffices to show that the dimensions of the semi-abelian varieties~$B$
such that there exists a maximal $k$\nobreakdash-morphism $b \colon X \rightarrow B$ mapping~$x$ to~$0$ are bounded.
Let $b \colon X \rightarrow B$ be such a morphism.  Let $u \colon X \otimes_k \kalg \rightarrow A$ be the universal
morphism given by Theorem~\ref{thalbexiste} over~$\kalg$ (so~$A$ is a semi-abelian variety over~$\kalg$).  Let $h \colon A \rightarrow B \otimes_k \kalg$
be the unique $\kalg$\nobreakdash-homomorphism such that $h \circ u = b \otimes_k \kalg$.
The $k$\nobreakdash-morphism~$b$ is generating since it is maximal.  Therefore, by Lemma~\ref{lemmagenerating}, the $\kalg$\nobreakdash-morphism $h\circ u$
is generating.  It follows that~$h$ is surjective, hence $\dim(B)$ is bounded by~$\dim(A)$.
\end{sketchofproof}

\bigskip
By a standard Galois descent argument, it follows from Theorem~\ref{thalbexiste} applied to a separable closure of~$k$
that the Albanese torsor~$\Alb^1_{X/k}$ and the Albanese variety~$\Alb^0_{X/k}$ exist for geometrically integral varieties~$X$ over an arbitrary field~$k$ (see~§\ref{sectiontheorder}
for the definitions of these objects, and \cite[Theorem~2.1]{dennis} or~\cite[Ch.~V, §4]{serregalg} for the Galois descent argument).

Let~$X$ be a smooth proper geometrically integral variety over~$k$ and $U \subseteq X$ be a dense open subset.
Under the assumption that~$k=\kalg$, Serre gave in~\cite[Théorème~1]{serremunivdiff}
an explicit description of the semi-abelian variety~$\Alb^0_{U/k}$ in terms of the Picard variety of~$X$ and the group of divisors on~$X$
which are algebraically equivalent to~$0$ and supported on~$X \setminus U$.  It turns out that the proof given in~\cite{serremunivdiff} applies verbatim when~$k$ is only separably closed,
thus yielding, after a slight reformulation
(which was not possible at the time~\cite{serremunivdiff} was written; indeed~\cite{serremunivdiff} predates the definition of $1$\nobreakdash-motives):

\bigskip
\begin{theorem}[ (Serre)]%
\label{thserreexpl}
Let~$X$ be a smooth proper geometrically integral variety over a field~$k$ and $U \subseteq X$ be a dense open subset.
Let~$\kbar$ be a separable closure of~$k$ and~$D$ denote the étale $k$\nobreakdash-group scheme defined
by $D(\kbar)=\Div^0_{X \setminus U}(X \otimes_k \kbar)$.
The semi-abelian variety~$\Alb^0_{U/k}$ is canonically dual to the $1$\nobreakdash-motive $[D \rightarrow P]$,
where $P=\Pic^0_{X/k,\mathrm{red}}$ is the Picard variety of~$X$
and $D \rightarrow P$ is the natural map.
\end{theorem}

\bigskip
\begin{corollary}
\label{corformation}
We keep the hypotheses of Theorem~\ref{thserreexpl}.  Let~$K/k$ be a separable extension (not necessarily algebraic).
Then the natural morphism $\Alb^n_{U \otimes_k K/K} \rightarrow \Alb^n_{U/k} \otimes_k K$
is an isomorphism for every $n \geq 0$.
\end{corollary}

\bigskip
The separability hypothesis in Corollary~\ref{corformation} is known to be superfluous in the classical situation where $U=X$ (see, \emph{e.g.}, \cite[Théorème~3.3]{tdte6}).
However it cannot be removed in general.  Indeed, let~$k$ be any imperfect field, let $X=\P^1_k$, let $P \in \A^1_k$ be a closed point whose residue field is purely
inseparable over~$k$,
let~$K$ be the residue field of~$P$, and let $U=\A^1_k \setminus \{P\}$; then,
according to Theorem~\ref{thserreexpl},
the natural morphism of tori $\Alb^0_{U \otimes_k K/K} \rightarrow \Alb^0_{U/k} \otimes_k K$ is an isogeny of degree~$[K:k]$
(more precisely, the map of character groups it induces is isomorphic to multiplication by~$[K:k]$ on~$\Z$).

\let\section=\myoldsection
\bibliographystyle{amsplain}
\bibliography{alb}
\end{document}